\documentclass[12pt]{article}
\usepackage{amsfonts}
\usepackage{latexsym}

\newcommand{\qed}{$\hfill\Box$}

\newcommand{\R}{\mathbb{R}}

\newcommand{\E}{\mathbb{E}}
\newcommand{\PP}{\mathbb{P}}

\newcommand{\al}{\alpha}
\newcommand{\be}{\beta}

\newcommand{\ga}{\gamma}
\newcommand{\ka}{\kappa}
\newcommand{\Ga}{\Gamma}
\newcommand{\si}{\sigma}

\newcommand{\ep}{\varepsilon}

\newcommand{\de}{\delta}
\newcommand{\te}{\theta}

\newcommand{\De}{\Delta}
\newcommand{\om}{\omega}
\newcommand{\Om}{\Omega}
\newcommand{\ze}{\zeta}

\newcommand{\n}{\mathcal N}

\newcommand{\f}{\mathcal F}

\def\ni{\noindent}

\newcommand{\st}{\sum_{i=1}^{[t/\De_n]}}

\newcommand{\proba}{(\Omega ,\f,(\f_t)_{t\geq0},\PP)}

\newcommand{\Wb}{\widetilde{b}}

\newcommand{\umu}{\underline{\mu}}
\newcommand{\unu}{\underline{\nu}}

\newcommand{\bee}{\begin{equation}}
\newcommand{\eee}{\end{equation}}
\newcommand{\bea}{\begin{eqnarray}}
\newcommand{\eea}{\end{eqnarray}}
\newcommand{\bean}{\begin{eqnarray*}}
\newcommand{\eean}{\end{eqnarray*}}

\newcommand{\oti}{\otimes}
\newcommand{\ti}{\times}
\newcommand{\fon}{\Omega\times\R_+\times\R^d\times\R^d}

\newtheorem{prop}{Proposition}[section]

\newtheorem{defi}[prop]{Definition}
\newtheorem{lem}[prop]{Lemma}
\newtheorem{theo}[prop]{Theorem}
\newtheorem{rem}[prop]{Remark}

\begin{document}
\title{Convergence of some random functionals of discretized semimartingales}

\author{Assane Diop
\thanks{Laboratoire de Probabilit\'es et Mod\`eles Al\'eatoires, Universit\'e Pierre et Marie Curie,
 4 Place Jussieu \newline
75252 Paris Cedex 05, France. E-mail: assane.diop.math@gmail.com}}


\maketitle

\abstract{In this paper, we study the asymptotic behavior of sums of functions of the increments of  a
given semimartingale, taken along a regular grid whose mesh goes to $0$. The function of the $i$th
increment may depend on the current time, and also on the past of the semimartingale before this
time. We study the convergence in probability of two types of such sums, and we also give associated
central limit theorems. This extends known results when the summands are a function depending
only on the increments, and this is motivated mainly by statistical applications.}

\vspace{.3cm}

\ni{\it Keywords}: Contrast functions, Power variation, Limit theorems, Semimartingale.

\section{Introduction}
In many practical situations, one observes a random process $X$ at
discrete times and one wants to deduce from these observations, some
properties on  $X$. Take for example the specific case of a
$1$-dimensional diffusion-type process  $X=X^\te$ depending
on a real-valued parameter $\te$, that is:
\bee\label{e11}
dX_s~=~\si(\te,s)\,dW_s~+~a(\te,s)\,ds,
\eee
where $\si$ and $a$ are
(known) predictable functions on $\Om\ti\R_+$, and where $W$ is a
Brownian motion. We observe the values of $X$  at times
$i\De,~i=0,1,\,2,\cdots,n\De$, and the aim is to estimate $\te$. There
are two cases: in the first one the observation window is
arbitrarily large. In the second case (which is our concern here),
the observation window is fixed, and so $\De=\De_n$ goes to $0$ and
$T=n\De_n$ is fixed.

Most known methods rely upon
minimizing some contrast functions, like minus the log-likelihood,
and those are typically expressed as ``functionals''
of the form:
\bee\label{e4}
\sum_{i=1}^{n} g_n\left(\si(\te',(i-1)\De_n),X^\te_{(i-1)\De_n}
,
X^\te_{i\De_n}-X^\te_{(i-1)\De_n}\,\right),
\eee
with $g_n$ possibly depending on $n$, see for example \cite{JG}.
In other words, the asymptotic behavior (convergence, and if
possible associated central limit theorems) of functionals like
(\ref{e4}) is very important. This is why,  for a function
$f:\Om\ti\R_+\ti\R^d\ti\R^d\to\R$ and  a $d$-dimensional
semimartingale $X$, we study the asymptotic behavior of the following
two sequences of process
\bee\label{e12} \left. \begin{array}{lll}
V^{n}(f,X)_t & = & \st f\left(\om,
(i-1)\De_n,X_{(i-1)\De_n}, X_{i\De_n}-X_{(i-1)\De_n}\right),\\[2mm]
V^{'n}(f,X)_t & = & \De_n \st f\left(\om,(i-1)\De_n,X_{(i-1)\De_n},
\frac{X_{i\De_n}-X_{(i-1)\De_n}}{\sqrt{\De_n}}\right),
\end{array} \right\}\eee
when $\De_n\to0$. So, providing some basic tools
for statistical problems is our main aim in
this paper, although we do not study any specific statistical problem.

Another motivation for studying functionals like (\ref{e12})
is that they appear naturally
in numerical approximations of stochastic differential equations like
the Euler scheme or more sophisticated discretization schemes.

Let us now  make two comments on the third argument of $f$ in
the processes   in (\ref{e12}), namely $X_{(i-1)\De_n}$:
\begin{enumerate}
\item The functionals (\ref{e12}) are not changed if we replace $f$
by $g(\om,t,x)=f(\om,t,X_{t-}(\om),x)$, so apparently one could dispense
with the dependency of $f$ upon its third argument. However, we will
need some H\"older continuity of $t\mapsto g(\om,t,x)$ which
is {\em not}\ satisfied by $g$ defined as just above: so it is more
convenient to single out the third argument.

\item One could replace $X_{(i-1)\De_n}$ by $Y_{(i-1)\De_n}$ for another
semimartingale $Y$, say \\$d'$-dimensional. But those apparently more
general functionals are like (\ref{e12}) with the
$(d+d')$-dimensional pair $Z=(Y,X)$ instead of $X$.
\end{enumerate}
\vspace{.2cm}

When $f(\om,s,z,x) \equiv f(x)$ ($f$ is ``deterministic''),  (\ref{e12})
becomes:
\bee\label{e4'}
\left.\begin{array}{lll} V^n(f,X) & = &
\st f\left(X_{i\De_n}-X_{(i-1)\De_n}\right),\\[2mm]
 V'^{n}(f,X) & = & \De_n\st
 f\left(\frac{ X_{i\De_n}-X_{(i-1)\De_n}}{\sqrt{\De_n}}\right).\end{array}
 \right\}
\eee
When further $f(x)=|x|^r$, the processes $V^n(f,X)$ are known as
the {\it realized power variations}, and of course
$V^{\prime n}(f,X)=\De_n^{1-r/2}\,V^n(f,X)$.

The convergence of power variations
is not new, see for example \cite{L}, an old paper by L\'epingle.
Recently they have been the object of a
large number of papers, due to their
applications in finance. Those applications are essentially the
estimation of the volatility and tests for the presence or absence of jumps.

An early paper is Barndorff-Nielsen and Shephard \cite{BS},
when $X$ is a continuous It\^o's semimartingale.
Afterwards,  many authors studied these type of processes:
Man\-ci\-ni \cite{MAN} studied  the case where $X$ is discontinuous
with L\'evy type jumps,
in \cite{JJ1} Jacod studied the general  case of a L\'evy process,
 Corcuera, Nualart and  Woerner in \cite{CNW} studied the case of a
fractional process, ..., the list is far from exhaustive. The results
appear in their most general form for a continuous semimartingale
 in \cite{BGJPN} and a discontinuous one in \cite{JJ}.

To give an idea of the expected results, let us mention that when
$X$ is a $1$-dimensional It\^o's semimartingale with diffusion coefficient
$\si$
and when $f$ is continuous and ``not too large near infinity'' (depending
on whether $X$ is continuous or not) we have
 $$V'^{n}(f,X)_t \stackrel{\PP}{\longrightarrow}~
 \int_0^t\rho_{\si_s}(f)\,ds,$$
(see for example \cite{BGJPN}),
where $\rho_x$ is the law of the normal variable ${\cal N}(0,x^2)$ and
$\rho_x(f)$ is the integral of $f$ with respect to $\rho_x$.

In \cite{BS} Barndorff-Nielsen and Shephard give a central limit
theorem for $V'^{n}(f,X)$, using a result of
Jacod and Protter about a central limit theorem (or: CLT) for
the Euler scheme for stochastic differential equations, see \cite{JP}.
This CLT has been generalized in many papers, like
\cite{BGJPN} when $X$ is continuous. If $X$ is discontinuous,
Jacod (in \cite{JJ}) gives a  CLT  when the Blumenthal-Getoor
index $p$ of $X$ is smaller than $1$, and no CLT is known when $p>1$.

Concerning $V^n(f,X)$,  in the uni-dimensional case,
Jacod extends some old results of L\'epingle in \cite{L}. In particular,
if $f(x)\sim|x|^r$ near the origin and is continuous and $X$ is
an arbitrary semimartingale, then
\bee\label{e13}
V^n(f,X)~\stackrel{\PP}{\longrightarrow}~ D(f,X),\eee with
$$D(f,X)_t=\left\{ \begin{array}{lll}  \sum_{s\le t} f(\De X_s)~\mbox{if}~ r>2, \mbox{ or if }
r\in (1,2)\mbox{ and } \langle X^c, X^c \rangle \equiv0,
 \\[2mm]\sum_{s\le t} f(\De X_s) +
\langle X^c, X^c \rangle_t,  ~\mbox{if}~  r=2,
\end{array} \right. $$
  where $\De X_s$ is the jump of $X$ at time $s$, and $X^c$
denotes
the continuous martingale part of $X$.
Moreover, Jacod gives a  central limit theorem for $V^n(f,X)$, first for
L\'evy processes in \cite{JJ1}, second for semimartingales in \cite{JJ}.

The difficulty of the extended setting in the present paper is due to
the fact that $f$  is not any more
deterministic and depends on all the variables $(\om,s,z,x)$, as we
have seen in the statistical problem. We want to know to which
extent  the earlier results remain valid in this setting,
and especially the CLTs. Our concern is to exhibit reasonably
general conditions on the test function $f$ which ensure that
the previously known results extend. Note also that for the CLT
concerning $V'^n(f,X)$, and contrary to the existing literature, we do
not always assume that $f(\om,t,z,x)$ is even in $x$, although most
applications concern the even case. The reader will also observe
that in some cases there are additional terms due to the parameter
$z$ in $f(\om,t,z,x)$.

The paper is organized as follows: in Sections 2 and 3 we state the Laws
of large numbers and the CLT respectively, and in Sections 4 and
5 we give the proofs.

\vspace{.3cm}
\ni{\bf Acknowledgement}

Thank you Professor Jean Jacod for  introducing  me to this research theme
and for your constant and so valuable assistance  during
 the writing of this paper.

\section{Laws of large numbers}\setcounter{equation}{0}
\subsection{General notation}

The basic process $X$ is a $d$-dimensional semimartingale on a fixed
filtered probability space $\proba$. We denote by $\De X_s
=X_s-X_{s-}$ the jump of $X$ at time $s$, and by $I$ the set
$$I~=~\Big\{r\geq0:~\sum_{s\le t} ||\De X_s||^r\,<\,\infty ~~\mbox{a.s
for all $t$}\Big\}.$$
Note that the set $I$ always contains the interval $[2,\infty)$.

The optional and predictable $\si$-fields on $\Om\ti\R_+$ are denoted
by ${\cal O}$ and ${\cal P}$, and if $g$ is a function on $\Om\ti\R_+\ti
\R^{l}$ we call it optional (resp. predictable) if it is
${\cal O}\oti{\cal R}^l$-measurable
(resp. ${\cal P}\oti{\cal R}^l$-measurable), where ${\cal R}^l$ is
the Borel $\si$-field on $\R^l$.

The function $f$ (unless otherwise stated) denotes  a function from
$\fon$ into $\R^q$, for some $q\geq1$ .
When $f(\om,t,z,x)$ admits partial derivatives in $z$ or $x$, we
denote by $\nabla_zf$ or $\nabla_xf$ the corresponding gradients.

If $M$ is a matrix, its transpose is $M^t$. The set
of all $p\ti q$ matrices is ${\cal M}(p,q)$, and
${\cal T}(p,q,r)$ is the set of all $p\ti q\ti r$-arrays.

For any  $\si\in{\cal M}(d,m)$ we denote by $\rho_{\si}$ the normal
law $\n(0,\si\si^t)$, and by  $\rho_{\si}(f(\om,s,z,.))$ the integral of the
function $x\mapsto f(\om,s,z,x)$ with respect to $\rho_{\si}$.

We denote by ${\cal B}$ the set of all functions $\phi:\R^d\to\R_+$
bounded on compact.

A sequence $(Z^n_t)$ of processes is said to converge u.c.p. (for:
uniformly on compact sets and in probability) to $Z_t$, and
written $Z^n\stackrel{u.c.p}{\rightarrow} Z$
or $Z^n_t\stackrel{u.c.p}{\rightarrow} Z_t$,  if
$\PP\left(\sup_{s\le t}||Z_s^n-Z_s||>\ep\right)$ $\to0$ for all $\ep,t>0$.

We write $Z^n\stackrel{{\cal L}-(s)}{\rightarrow} Z$ or
$Z^n_t\stackrel{{\cal L}-(s)}{\rightarrow} Z_t$, if the process $Z^n$
converge stably in law to $Z$, as processes
(see \cite{JS} for details on the stable convergence).

We gather some important properties of $f$ in the following
definition.

\begin{defi}\label{d1}
a)  We  say that $f$ is of
(random) polynomial growth if there exist a locally bounded process $\Ga$
(meaning: $\sup_{s\leq T_n}\Ga_s\leq n$ for a sequence $T_n$ of
stopping times increasing a.s. to $\infty$), a
function $\phi\in {\cal B}$, and a real $p\geq0$ such that
\bee\label{er1}||f(\om,s,z,x)||\le
\Ga_s(\om)\phi(z)(1+||x||^p).\eee
If we want to specify $p$, we say that $f$ is at most of $p$-polynomial growth.

b) we  say that $f$ is locally equicontinuous in $x$ (resp.
$(z,x)$)  if for all $\om$, all $T>0$, and all compacts ${\cal K,K'}$
in $\R^d$,
the family of functions $\big(x\mapsto f(\om,s,z,x))_{s\leq T,z\in
{\cal K'}}$ (resp. $\big((z,x)\mapsto f(\om,s,z,x))_{s\leq T}$)
is equicontinuous on ${\cal K}$ (resp. ${\cal K}\ti{\cal K'}$).
\end{defi}

\subsection{Assumptions}\label{s4}
Let us start with the assumptions on $X$. For $V^n(f,X)$  we only
need $X$ to be an arbitrary semimartingale. For $V'^n(f)$ we need
$X$ to be an It\^o semimartingale and a little more.
Recall first that the property of $X$ to be an It\^o semimartingale
is equivalent to the following: there are, possibly on an extension of the
original probability
space, an $m$-dimensional Brownian motion $W$ (we may always take $m=d$)
and a Poisson random measure $\umu$ on $\R_+\ti\R$ with intensity
measure $\unu(ds,dy)=F(dy)\,ds$  with $F$ is a $\si$-finite measure on
$\R$, such that $X$ can be written as
\begin{eqnarray}
X_t  & = &  X_0 +\int_0^t b_s\,ds +\int_0^t\si_{s-}
dW_s+\int_0^t\int_{\R}h\left(\de(s,y)\right)(\umu-\unu)(ds,dy) \nonumber\\
 & & +  \int_0^t\int_{\R} h'\left(\de(s,y)\right)\umu(ds,dy),\label{e5}
\end{eqnarray}
for suitable ''coefficients'' $b$ (predictable $d$-dimensional),
$\si$ (optional $d\ti m$-dimensional), $\de$ (predictable $d$-dimensional
function on $\Om\ti\R_+\ti\R$) and $h$ is a truncation function
from $\R^d$ into itself (continuous with compact support, equal to
the identity on a neighborhood of $0$), and $h'(x):=x-h(x)$.

Then we set:

\vspace{0.5cm}
\ni{\bf Hypothesis $(N_0)$:}
The process $X$ is an It\^o's semimartingale, and its coefficients
in (\ref{e5}) satisfy the following:
$b$ and $\int_{\R}(1\land||\de(\om,s,y)||^2)\,F(dy)$
are locally bounded, and $\si$ is c\`adl\`ag.\qed
\vspace{0.5cm}

For the test function $f$ we introduce the following, where
$A$ is an arbitrary subset of $\R^d$:

\vspace{.5cm}

\ni {\bf Hypothesis $(K[A])$:} $f(\om,t,z,x)$ is continuous in
$(z,x)$ on $\R^d\ti A$ and if $(t_n,z_n,x_n)\to (t,z,x)$ with $x\in A$ and
 $t_n<t$, then $f(\om,t_n,z_n,x_n)$ converges to a limit depending
on $(\om,t,z,x)$ only, and denoted by $f(\om,t-,z,x)$. \qed

\vspace{.5cm}

\subsection{Results}
The  first two theorems concern the processes $V^n(f)$.

\begin{theo}\label{t1} Let $X$ be an arbitrary semimartingale, and let $f$
satisfy $K(\R^d)$. Suppose  there exist
a neighborhood $V$ of $0$ on $\R^d$, a real $p>2$, and for any  $K>0$, a
locally bounded process  $\Ga^K$ such that:
\bee\label{e6}
\|z\|\le K,~\,x\in V\,\Rightarrow\,
||f(\om,s,z,x)||~\le~\Ga^K_s(\om)\|x\|^p.
\eee
Then $V^n(f)$ converge a.s.\ for the Skorokhod topology to the process
\bee\label{e16} D(f)_t~=~\sum_{s\le t} f(s-,X_{s-},\De X_s).\eee
\end{theo}

\begin{rem}\label{r5}
This is one of the rare situations where one has  almost sure convergence;
see Section 3.1 of  \cite{AD} for some  other ones.
\end{rem}

\begin{theo}\label{t2} Let $X$ be an arbitrary semimartingale, and let
$f$ be optional, satisfy $(K(\R^d))$ and
$f(\om,s,z,0)=0$, and be $C^2$ in $x$ on some neighborhood $V$ of
$0$, and assume also
\begin{itemize}
\item[$\bullet$]For any $j,k \in \{1,\cdots,d\}$, the functions
$\frac{\partial f}{\partial x_j}(\om,s,x,z)$ and
$\frac{\partial^2 f}{\partial x_j
\partial x_k}(\om,s,x,z)$ defined on $\Om\ti\R_+\ti\R^d\ti V$
satisfy $(K[V])$.
\item[$\bullet$] There exist $\phi \in{\cal B }$ and a locally bounded
process $\Ga$   such that
 $$\sum_{j=1}^d\left(\left\|\frac{\partial f}{\partial
x_j}(s,z,0)\right\|~+
 ~\sum_{k=1}^d\,\left(\sup_{x\in\,V}\,\left\|
\frac{\partial^2f}{\partial x_j \partial x_k} (s,z,x) \right\|\right)
\right)~\le~\Ga_s\phi(z).
$$
\end{itemize}
Then $V^n(f)$ converge in probability, in the Skorokhod sense, to the
process
\begin{eqnarray}
D(f)_t &= & \sum_{j=1}^d\int_0^t\frac{\partial f}{\partial
x_j}(s-,X_{s-},0)\,dX_s +\frac1{2}\sum_{j,k=1}^d
\int_0^t\frac{\partial^2f}{\partial
x_j \partial x_k}(s-,X_{s-},0)\,d\langle X^{j,c},X^{k,c}
\rangle_s\nonumber\\[2mm]
& & +\sum_{0<s\le t}\left(f(s-,X_{s-},\De X_s)- \sum_{j=1}^d
\De X_s^j\frac{\partial
f}{\partial
x_j}(s-,X_{s-},0)\right),\label{e15} \end{eqnarray}
where $X^c$ is the continuous martingale part of $X$.
\end{theo}

The two versions (\ref{e16}) and (\ref{e15}) of $D(f)$ agree when
$f$ satisfies the hypotheses of Theorem \ref{t1}, so Theorem \ref{t2}
extends Theorem \ref{t1} and gives the results in a more
complete form. This result was not known even in the case when $f$
only depends on $x$.

\begin{rem}\label{r4}
Both theorems remain valid if the
discretization grid is not regular, provided the successive discretization
times are stopping times and the mesh goes to $0$
(see  Sections 3.5 and 4.5 of \cite{AD} for results of this type).
\end{rem}

Now we state the result about $V'^{n}(f)$.

\begin{theo}\label{t3}
Let  $f$ be optional, satisfy $(K(\R^d))$, be locally equicontinuous in $x$
and with $p$-polynomial growth. Assume further that one of the following two
conditions is satisfied:
\begin{enumerate}
\item $X$ satisfies $(N_0)$ and  $p<2$.
\item  $X$ satisfies $(N_0)$ and is  continuous.
\end{enumerate}
Then \bee\label{ep30}{V'}^n(f)~\stackrel{u.c.p.}{\longrightarrow}~
\int_0^t\,\rho_{\si_{s-}}\left(f(s-,X_{s-},.)\right)\,ds.\eee
\end{theo}

\begin{rem}\label{r3}
Comparing with \cite{BGJPN} or  \cite{JJ}, we see that there is no
additional term due to the third
argument $z$ in $f(\om,s,z,x)$.

In the discontinuous case (Hypothesis 1), the condition $p<2$  simplifies
the computations but is not optimal. The result remains true
valid if there exist   $\phi,\,\phi'\,\in\,{\cal B}$ such that:
$$\phi'(x)\to 0,~\mbox{when}~||x||\to\infty,~
\mbox{and}~||f(\om,s,z,x)||~\le~\Ga_s(\om)\phi(z)||x||^2\phi'(x).$$
\end{rem}

\section{Central limit theorems}\label{ss1}

In the framework of the CLT, one needs some additional assumptions
both on $X$ and on $f$, which depend on the
problem at hand.

\subsection{Assumptions on $X$}

\ni{\bf Hypothesis $(N_1)$:}
$(N_0)$ is satisfied, and there exist a sequence $(S_k)$ of stopping times
increasing to $\infty$ and deterministic Borel functions $(\ga_k)$ such that:
$$\left|\left|\de(\om,s,y)\right|\right|~\le~\ga_k(y)~\mbox{if}~s\le
S_k(\om)~~~\mbox{and}~~\int_{\R} (1\land\ga_k(y)^2)\,F(dy)~<\infty.$$
\vskip-9mm \qed
\vspace{0.5cm}

The next assumption depends on a real $s\in[0,2]$:

\vspace{0.5cm}

\ni{\bf Hypothesis $(N_2(s))$:}
$(N_1)$ is satisfied, the mapping $s\mapsto\de(\om,s,y)$ is c\`agl\`ad, and
$\int_{\R} (1\land\ga_k(y)^s)\,F(dy)~<\infty$.
Moreover, the process $\si$ in (\ref{e5}) satisfies:
\bee\label{e14} \si_t ~= ~\si_0 + \int_0^t\widetilde{b}_u
du+\int_0^t\widetilde{\si}_u
dW_u+ M_t+ \sum_{u\le t} \De \si_u 1_{\{||\De \si_u||\ge 1\}},\eee
where
 \begin{itemize}
    \item[$\bullet$]  $\widetilde{b}$ is predictable and locally
    bounded.
\item[$\bullet$] $\widetilde{\si}$ is  c\`adl\`ag, adapted with
values in ${\cal T}(d,m,m)$.

\item[$\bullet$]  $M$ is an ${\cal M}(d,m)$-valued local martingale,
orthogonal to $W$  and satisfying $||\De M_t||\le 1$ for all $t$.
Its predictable quadratic covariation is
$\langle M,M\rangle_t=\int_0^t a_u\,du,$, where $a$ is
locally bounded.

\item[$\bullet$] The predictable compensator of
$\sum_{u\le t}  1_{\{||\De \si_u||\ge 1\}}$ is
$\int_0^t \widetilde{a}_u\,du$, where $\widetilde{a}$ is locally bounded.\qed
\end{itemize}
\vskip4mm

Clearly  $(N_2(s))\,\Rightarrow\,(N_2(s'))$, if $s<s'$.

\begin{rem}\label{r1} It is well known that the assumptions on $\si$
in $(N_2(s))$ may be replaced by the following one (up to modifying the
Poisson measure $\umu$):
\begin{eqnarray}\si_t & = & \si_0+\int_0^t\widetilde{b}_u\, ds+\int_0^t
\widetilde{\si}_u
dW_u+\int_0^t\widetilde{v}_u\,dV_s+\int_{\R}\int_0^t k(\widetilde{\de}(u,y))
\star(\umu-\unu)(du,dy) \nonumber \\[2mm]
  & & + \int_{\R}\int_0^t k'(\widetilde{\de}(u,y))\star\umu(du,dy),\label{e9}
  \end{eqnarray}
where $\widetilde{b}$ and $\widetilde{\si}$ are like in $(N_2(s))$ and
\begin{itemize}
\item $V$ is a $l$-dimensional
Brownian motion independent of $W$.
\item $\widetilde{v}$ takes its values in ${\cal T}(d,m,l)$,
is progressively measurable and locally bounded.
\item $k(x)$ is a  truncation function on $\R^d\ti\R^m$ and $k'(x):= x-k(x)$.
\item $\widetilde{\de}:\,\Om\ti\R_+\ti\R\to{\cal M}(d,m)$ is predictable and is such that:
$\int_{\R} (1\land ||\widetilde{\de}(u,y)||^2)\,F(dy)$ is locally bounded.
\end{itemize}
Of course, $a,~\widetilde{a},~\widetilde{v}$ and $\widetilde{\de}$ are
related, for example  if $k(x)= x 1_{\{||x|| < 1\}}$, one has
$\widetilde{v}^2_u+\int_{\{||\widetilde{\de}(u,y)||\leq1\}}
 {\widetilde{\de}}^2 (u,y)\,F(dy)~=~a^2_u$
 and $\widetilde{a}_u=\int_{\{||\widetilde{\de}(u,y)||>1\}} \,F(dy)$.
\end{rem}

\subsection{Assumptions on the test function $f$}

\ni{\bf Hypothesis ($M_1$):}
$f$ is optional and there exists a neighborhood $V$ of $0$ such that
$f(\om,s,z,x)$ is $C^1$ in $(z,x)$, the functions   $\nabla_x f$,
$\nabla_z f$ are $C^1$ in $x$ on $ V$, and
$$f(\om,s,z,0)=\nabla_x f (\om,s,z,0)\equiv 0.$$
Moreover there are  a locally
bounded process $\Ga$, a real $\al>\frac1{2}$, and
some  functions $\phi,\,\ep$ and $\te$ belonging to ${\cal B}$,
with $\ep(x)\to 0$ as $||x||\to 0$  and $\te(x)\leq\|x\|^2$ in the
neighborhood of $0$,  such that:
$$ \sum_{j,j'=1}^d
  \left(\left\| \frac{\partial^2 f}{\partial x_j
  \partial x_{j'}}(\om,s,z,x)\right\|
+\left\| \frac{\partial^2 f}{\partial x_j
\partial z_{j'}}(\om,s,z,x)\right\|\right)\,
\le\Ga_s(\om)\phi(z)||x||\ep(x),$$
\ni and for all $T>0$ and $s,t\in[0,T]$,
\bee\label{e7}
\left\|f(\om,t,z,x)-f(\om,s,z,x)\right\| \le
\Ga_T(\om)\phi(z)\, |t-s|^{\al}\;\te(x).
\eee
\vskip-3mm\qed

\vspace{.5cm}

\ni{\bf Hypothesis $(M_2)$:}
$f(\om,t,z,x)$ is optional, $C^1$ in $(z,x)$, with $\nabla_xf$
and $\nabla_zf$ of (random) polynomial growth and locally equicontinuous
in $(z,x)$, and there are $\Ga$, $\phi$, $\al$ as in $(M_1)$
and some $p>0$ such that for all $T>0$ and $s,t\in[0,T]$,
\bee\label{ep23}
\|f(\om,s,z,x)-f(\om,t,z,x)\|~\le~
 \Ga_T(\om)\phi(z)|t-s|^{\al}(1+\|x\|^p),
\eee
\vskip-3mm\qed

\vspace{.5cm}

\ni{\bf Hypothesis $(M_2')$:}
 $(M_2)$ is satisfied  and moreover
 $$\|f(\om,s,z,x)\|+\|\nabla_x f(\om,s,z,x)\|~\le~\phi(z)\Ga_s(\om).$$
\vskip-8mm\qed

\vspace{.5cm}

The previous hypotheses are fulfilled by most of the test
functions used in statistics.

\subsection{The results}

In order to define
the limiting processes, we need to expand the original space $\proba$,
 what we do as follows:

Consider an auxiliary space $(\Om',\f',\PP')$, which supports a
$q$-dimensional Brownian motion $\overline{W}$ and some sequences
$\big\{(U^k_p)_{1\le k\le m};~(U^{'k}_p)_{1\le k\le m};~
(\kappa_p)\big\}_{p\ge 1}$ of random variables, where the $U^k_p$ and
$U^{'k}_p$ are normal   ${\cal N}(0,1)$  and the   $(\kappa_p)$ are
uniform  on $(0,1)$. We suppose  all these variables and processes mutually
independent.

Now set:
$$ \widetilde{\Om}=\Om\ti\Om',~~\widetilde{\f}=\f\oti\f',~~
\widetilde{\PP}=\PP\otimes\PP'.$$
We then extend the variables and processes defined on $\Om$ or
$\Om'$ on the space $\widetilde{\Om}$, in the usual way.

 Let $(T_p)$ be an arbitrary sequence of stopping times exhausting the
jumps of $X$ (meaning: they are stopping times such that
for all $(\om,s)$ with $\De X_s(\om)\ne 0$, there exists a
unique $p$ such that $T_p(\om)=s$).
We define on $\widetilde{\Om}$ the filtration $(\widetilde{\f}_t)$ which is the
smallest one satisfying the following conditions:
\begin{itemize}
\item  $(\widetilde{\f}_t)$ is right continuous, and
$\f_t\subset \widetilde{\f}_t$,
\item  $\overline{W}$ is adapted on $(\widetilde{\f}_t)$,
\item the variables $U^{k}_p,~U^{'k}_p$ and $\kappa_p$
are  $\widetilde{\f}_{T_p}$ measurable.
\end{itemize}

Now we are ready to give the results. We start with $V^n(f)$:

\begin{theo}\label{t4}
Suppose that $X$ satisfies $(N_1)$ and  $f$ satisfies $(M_1)$,
then
$$\frac1{\sqrt{\De_n}}\left(V^n(f)-D(f)_{[t/\De_n]\De_n}\right)~
\stackrel{{\cal L}}{\longrightarrow}~F_t,$$
where the process $F$ is
\begin{eqnarray} F_t & = & \sum_{p:~T_p\le t} \sum_{j=1}^d
\sum_{k=1}^m\left( \left(
\sqrt{\ka_p}\, \si^{j,k}_{T_p-}U_p^k+\sqrt{1-\ka_p}\,
\si^{j,k}_{T_p}{U}_p^{'k}\right)
\frac{\partial f}{\partial x_j}\left(T_p-,X_{T_p-},\De X_{T_p}\right)
\right. \nonumber\\[2mm]
& & \left. - \sqrt{\ka_p}\,\si^{j,k}_{T_p-}U_p^k\frac{\partial f}
{\partial z_j}\left(T_p-,X_{T_p-},\De X_{T_p}\right)\right).
\label{ep15}\end{eqnarray}
\end{theo}

\begin{rem}\label{r2}
The last term in (\ref{ep15}) is due to the third argument of $f$, and
does not appear in \cite{JJ}. One could show that the theorem remains
valid  if, in the formula (\ref{e7}), $\te(x)\leq\|x\|^p$ near the origin
for some  $p\in[0,2]\cap I$.
\end{rem}

It is useful to give some properties of the process $F$ above. For this,
under $(M_1)$ and $(N_1)$, one defines an ${\cal M}(q,q)$-valued
process $C(f)$ as follows:
\begin{eqnarray} C(f)_t & = & \frac1{2}\sum_{p:~T_p\le t}\,\sum_{j,j'=1}^d\,
\sum_{k=1}^m\left\{
\left(
\si_{T_p-}^{j,k}\si_{T_p-}^{j',k}+\si_{T_p}^{j,k}\si_{T_p}^{j',k}\right)
\right.\nonumber \\[2mm]
 & & \ti \left(\frac{\partial f}{\partial x_j}\right)\left(\frac{\partial f}
 {\partial x_{j'}}\right)^t\circ(T_p-,X_{T_p-},\De X_{T_p})\nonumber\\[2mm]
 & & -\si_{T_p-}^{j,k}\si_{T_p-}^{j',k}\left(\left(\frac{\partial f}
 {\partial x_j}\right)
 \left(\frac{\partial f}{\partial z_{j'}}\right)^t
 +\left(\frac{\partial f}{\partial z_j}\right)
 \left(\frac{\partial f}{\partial x_{j'}}\right)^t\right)
\circ(T_p-,X_{T_p-},\De X_{T_p}) \nonumber\\[2mm]
& & \left.+ \si_{T_p-}^{j,k}\si_{T_p-}^{j',k}\left(\frac{\partial f}
{\partial z_j}\right)
\left(\frac{\partial f}{\partial z_{j'}}\right)^t
\circ(T_p-,X_{T_p-},\De X_{T_p}) \right\},
\label{e21}\end{eqnarray}

The following  lemma is given without proof, since it is an immediate
generalization of  lemma 5.10 of \cite{JJ}.

\begin{lem}\label{l3}
If  $(M_1)$ and  $(N_1)$ are satisfied, then $C(f)$ is well defined and
$F$ is a semimartingale on the extended space
$(\widetilde{\Om},\widetilde{\f},\widetilde{\PP})$.
If further   $C(f)$ is locally integrable, then $F$ is a locally
square-integrable martingale.

Conditionally on $\f$, the process $F$ is a square integrable
centered martingale with independent increments,
its conditional variance is
$C(f)_t=\widetilde{\E}\{F_t^2|\f\}$, its   law is completely
characterized by $X$ and $\si\si^t$ and  does not depend on the choice of
the sequence $(T_p)$.
\end{lem}

Now we turn to $V'^n(f)$.
Under $(M_2)$ or $(M_3(r))$, one defines a process $a$
taking its value in ${\cal M}(q,q)$ and satisfying for any
$j,k\in\{1,\cdots,q\}$:
\bee\label{e17}
\sum_{l=1}^q\,a^{j,l}_t\,a^{l,k}_t\,=\,\rho_{\si_{t}}
\left((f^jf^k)(t,X_t,.)\right)
-\rho_{\si_{t}}\left(f^j(t,X_t,.)\right)\rho_{\si_{t}}
\left(f^k(t,X_t,.)\right).
\eee
The  process $a$, which may be chosen $(\f_t)$-adapted, is the square-root of
the symmetric semi-definite positive element of
${\cal M}(m,m)$ whose components are given by the right side of (\ref{e17}).

\begin{theo}\label{t5}
Suppose   $f(\om,s,z,x)$ even in $x$, and assume  that one of the following
hypothesis is satisfied:
\begin{itemize}
\item $X$ is continuous and satisfies $(N_2(2))$  and $f$ satisfies $(M_2)$.
\item one has  $(N_2(s))$ for some $s\le 1$ and  $(M_2')$.
\end{itemize}
Then  $$\frac1{\sqrt{\De_n}}\left(V'^{n}(f)_t~-~
\int_0^t\,\rho_{\si_{s}}f(s,X_{s},.)\,ds~
\right)\,\stackrel{{\cal L}-(s)}{\longrightarrow}\, L(f)_t, $$
where
\bee\label{e35}
L(f)_t~=~\int_0^t a_s\,d\overline{W}_s.
\eee
\end{theo}

\begin{rem}
Some times, one wants to apply the theorem for functions of the type
$f(\om,s,z,x)=g(\om,s,z)\|x\|^r$,
which are not any more $C^1$ in $x$ on $\R^d$ when $r\in(0,1]$. Specifically,
consider  the following hypothesis:

\vspace{0.5cm}

\ni{\bf Hypothesis ($M_3(r)$):} $f(\om,s,z,x)$ is optional and there is a closed subset $B$ of $\R^d$
with Lebesgue measure $0$ such that the application $x\to f(\om,t,z,x)$ is
$C^1$ on $B^c$. Moreover there are $p\ge 0$ and $\al$,
$\phi$ and $\Ga$ as in $(M_1)$  such that for all $T>0$ and $s,t\in[0,T]$,
\bee\label{e1'}
\left. \begin{array}{lll} \|f(\om,s,z,x_1+x_2)-f(\om,s,z,x_1)\| & \le  &
\Ga_T(\om)\phi(z)\left(1+\|x_1\|^p\right)\|x_2\|^r.\\[2mm]
\|f\left(\om,s,z,x\right)-f\left(\om,t,z,x\right)\| & \le &
 \Ga_T(\om)\phi(z)|t-s|^{\al}\left(1+||x||^p\right). \end{array} \right\} \eee
 Moreover,
\begin{itemize}
\item[$\bullet$] if $r=1$ then $\nabla_x f$ defined on
$\Om\ti\R_+\ti\R^d\ti B^c$ is locally equicontinuous in $(z,x)$ with at most
polynomial growth.
\item[$\bullet$] if $r\ne 1$, then for any element $C\in{\cal M}(d,d)$
and any ${\cal N}(0,C)$-random vector $U$, the distance from $U$ to $B$
has a density $\psi_C$ on $\R_+$, satisfying
 $\sup_{x\in\R_+,\,\|C\|+\|C^{-1}\|< K}\,\psi_{C}(x)\,<\,\infty$ for all $K<\infty$.
 For any $x_1\in B^c,$
\bee\label{e2} \|\nabla_{x} f(\om,s,z,x_1)\| ~\le~\frac{\Ga_s(\om)\phi(z)(1+\|x_1\|^p)}{d(x_1,B)^{1-r}},\eee

and if  $||x_2||<\frac{d(x_1,B)}2$, then
      \bee\label{e3}\|\nabla_x f(\om,s,z,x_1+x_2)~-~\nabla_x f(\om,s,z,x_1)\| ~\le~\frac{\Ga_s(\om)\phi(z)(1+\|x_1\|^p)\|x_2\|}{d(x_1,B)^{2-r}}.\eee
     \end{itemize}
\vskip-2mm\qed

\noindent Then one can show that the results of theorem \ref{t5} remain valid
if $f$ satisfies $(M_3(r))$ for some $r\in(0,1]$ and
$X$ satisfies $N_2(2)$ with $\si\si^t$ everywhere invertible,
if further one of the following condition is satisfied:
\begin{itemize}
\item $f$ satisfies $(M_3(r))$  and $X$ is continuous,

\item $f$ satisfies $(M_3(r))$ and the real $p$ in (\ref{e1'}), (\ref{e2}) and
(\ref{e3}) is always equal $0$, while $X$ satisfies $(N_2(s))$ and either
$s\in [0, \frac2{3})$ and $r\in (0, 1)$ or
 $s\in (\frac2{3},1)$ and $r\in (\frac{1-\sqrt{3s^2-8s+5}}{2-s},1)$.
\end{itemize}
\end{rem}

\vspace{.5cm}

Our next objective is to generalize the CLT for $V'^{n}(f)$ in the case where
$f$ is not even. For this, we need some additional  notation.

Let  $U$ be  an  ${\cal N}(0,Id_m)$ random vector,  where $Id_m$ is the
identity matrix of order $m$ (recall that $m$ is  the dimension of the
Brownian motion $W$ in  $(N_2(s))$).
We then denote  by $\rho'$, the law of $U$ and by $\rho'(g_1(.))$ the
integral of any function
$g_1:\R^m\to\R^q$ with respect to $\rho'$ if it exists. If now
$g_2\,:\R^d\to\R^q$ and  $x\in  {\cal M}(d, m)$,
we  set:  $\rho'(g_2(x.))=\E\{g_2(xU)\}$.

For any $j\in \{1,\cdots,m\}$, we define the  projection $P_j$ on $\R^m$ by:
 $$P_j(u):=u_j~~\mbox{if}~~u=(u_1,\cdots,u_m).$$
Under $(M_2)$ we define  $w(1)$ and $w(2)$,  two adapted processes taking
their
values respectively in the spaces ${\cal M}(q, m)$ and  ${\cal M}(q, q)$,
and such that for all $j,k\in\{1,\cdots,q\}$ and $j'\in\{1,\cdots,m\}$
we have
\bee\label{ep31} \left.\begin{array}{l} w(1)^{j,j'}_s ~ =~
\rho'\left(f^j(s,X_s,\si_{s}.)P_{j'}(.)\right),\\[2mm]
\sum_{l=1}^q w(2)^{j,l}_t\, w(2)^{l,k}_t ~ =~ \rho'
\left((f^jf^k)(s,X_s,\si_s .)\right)\\\quad
- \rho'\left(f^j(s,X_s,\si_s .)\right)
\,\rho'\left(f^k(s,X_s,\si_s .)\right)-\sum_{l'=1}^m w(1)^{j,l'}_t\,
w(1)^{l',k}_t. \end{array} \right\}\eee
The process $w(2)$ is the square-root of the matrix whose components are
given by the right side of the second equality in (\ref{ep31}).
Finally, under $(N_2(2))$ set
\bee\label{e10} b'~=~b\,-\,\int_{\R}h(\de(s,y))\,F(dy).
\eee

\begin{theo}\label{t6}
Assume either one of the following two assumptions:
\begin{itemize}
\item $X$ satisfies $(N_2(2))$ and is continuous and $f$ satisfies $(M_2)$.
\item We have $(N_2(s))$ for some  $s\le 1$ and $f$ satisfies $(M_2')$.
\end{itemize}
If further $b'\equiv0$ and $\widetilde{\si}\equiv0$, we have
$$\frac1{\sqrt{\De_n}}\left(V^{'n}(f)_t\,-\,\int_0^t\,
\rho'f(s,X_{s},\si_{s}.)\,ds\,
\right)\,\stackrel{{\cal L}-(s)}{\longrightarrow}\, L(f)_t, $$
where
\bee\label{e20}
L(f)_t~:=~\int_0^t w(1)_s\,dW_s+\int_0^t w(2)_s\,d\overline{W}_s.
\eee
\end{theo}

\begin{rem}
Clearly, when $f$ is even in $x$, the two versions of the process $L(f)$
in Theorems \ref{t5} and \ref{t6}, agree.
If $X$ satisfies $(N_2(s))$ with $s\le 1$, the hypotheses
$b'=0$ and $\widetilde{\si}=0$ yield that $X$ has the form:
\bee\label{15} X_t= X_0+\int_0^t\si_s\,dW_s+\sum_{s\le t} \De X_s. \eee
\end{rem}

\section{Proof of the laws of large numbers}\setcounter{equation}{0}
\subsection{Theorems \ref{t1} and \ref{t2}}

We  start by stating two important lemmas, without proof. The first one
is a (trivial) extension of what is done in Subsection 3.1 of \cite{JJ},
and the hypothesis $(K(\R))$ plays a crucial role there. The second one
is a generalization of It\^o's formula, and its proof
can be found for example  in \cite{AD} (see lemma 3.4.2).

\begin{lem}\label{l1}
Let $X$ be an arbitrary semimartingale, and $f$ be a function satisfying
$(K[\R])$ and such that
$f(s,z,x)=0$ if $||x||\le \ep$ for some $\ep>0$.
Then $$V^n(f)_t~-~\sum_{s\le[t/\De_n]\De_n}~f(s-,X_{s-},\De X_s).$$
  converges in variation to $0$ when $n\to\infty$, for each $\om\in\Om$.
\end{lem}

\begin{lem}\label{l2}
Let $X$ be a semimartingale and $f(\om,u,z,x)$ be an optional function, $C^2$
in  $x$. Then for any $u$,
for almost all  $\om$ and for any $ t\ge u$, one has:
\begin{eqnarray*} f(u,X_u,X_t) & = & f(u,X_u,X_u) +
\sum_{j=1}^d\int_{u+}^t\frac{\partial f}{\partial x_j}(u,X_u,X_{s-})\,dX_s
\nonumber\\[2mm]
& & +\sum_{j,j'=1}^d\frac1{2}\int_{u+}^t\frac{\partial^2 f}
{\partial x_j\partial x_{j'}}(u,X_u,X_{s-})\,d\langle
X^{j,c},X^{j',c}\rangle_s\nonumber\\[2mm]
& & +  \sum_{u<s\le t}\Big(f(u,X_u,X_s)-f(u,X_u,X_{s-})-
\sum_{j=1}^d\De X^j_s \frac{\partial
f}{\partial x_j}(u,X_u,X_{s-})\Big).
\end{eqnarray*}
\end{lem}

Now we are ready to prove the two theorems about $V^n(f)$.
\vskip5mm

\noindent{\bf Proof of Theorem \ref{t1}:}
Since for any c\`adl\`ag process $Y$, the processes $Y_{[t\De_n]\De_n}$
converge pathwise to $Y$ for
the Skorokhod topology, it is sufficient  to prove that the processes
$V^n(f)_t-D(f)_{[t/\De_n]\De_n}$ converge u.c.p. to $0$.

We suppose  first that $\|X_t\|\leq C$ identically for some constant $C$.
Let  $t>0$, and $S^n=\{0=t_1^n<t_2^n<\cdots<t^n_{k^n}=t\}$ be a sequence
of partitions of
$[0,t]$ such that $\sup_{i}\, |t_i^n-t_{i-1}^n|~\to 0,$ when $n\to\infty$.
According to Th\'eor\`eme 4 of \cite{L}, one has:
$$\sum_{i=1}^{k^n}\,|X^j_{t_i^n}-X^j_{t_{i-1}^n}|^p~\longrightarrow~
\sum_{s\le t} |\De X_s^j|^p~~\mbox{a.s.},$$
for any $j\in\{1,\cdots,d\}$, and where $X^j$ is the $j$th component of $X$.

Since the mappings $t\mapsto\sum_{s\le t} |\De X^j_s|$
and $t\mapsto\st|\De_i^n X^j|^p$ are increasing, we deduce that
for almost all $\om$ and for any  real $t>0$,
\begin{equation}
\limsup_n \st \|\De_i^n X\|^p ~ \le~d^{p-1}
\sum_{j=1}^d \sum_{s\le t} |\De X_s^j|^p.\label{ep36}
\end{equation}

Let now  $\psi:~\R\to\R$ be a $C^{\infty}$ function such that
$1_{[-1,1]}(y)\le\psi(y)\le 1_{[-2,2]}(y)$. We then put for $y\in\R$
and $x\in\R^d$ and $\ep>0$:
\bee\label{e8} \psi_{\ep}(y)=
\left\{\begin{array}{ll}\psi(\frac{y}{\ep})& \mbox{if}~\ep<\infty\\
1&\mbox{if}~\ep=\infty,
\end{array} \right.\qquad
 \Psi_{\ep}(x)=\Pi_{j=1}^d \psi_{\ep}(x_j).
 \eee
Note that $$ \Psi_{\ep}(x)=
\left\{\begin{array}{ll}1 & \mbox{if}~||x||\le \ep\\
0&\mbox{if}~||x||>2d\ep,\end{array} \right. $$
and set, with the notation (\ref{e16}):
\bee\label{ep37}
Z^n(f)_t~=~V^n(f)_t-D(f)_t,
\eee
Then
\bee\label{ep38}
Z^n(f)=Z^n(f\Psi_{\ep})+Z^n(f(1-\Psi_{\ep})),\eee
\bee\label{ep39}
\limsup \sup_{t\le T} \|Z^n(f)_t\|  \le  \limsup
\sup_{t\le T} \|Z^n(f\Psi_{\ep})_t\|\, +\, \limsup \sup_{t\le T}
\|Z^n(f(1-\Psi_{\ep}))_t\|,\eee
for any $T>0$. By Lemma \ref{l1}, one has
\bee\label{ep40}\lim_{\ep\to 0} \, \limsup_n \, \sup_{t\le
T}~ \|Z^n(f(1-\Psi_{\ep}))_t\|~=~0.\eee
On the other hand, if $q\in(2,p)$ we have by (\ref{e6})
and $\|X\|\leq C$ and (\ref{ep36}):
\begin{eqnarray*}
\sup_{s\le t}~\|Z^n(f\Psi_{\ep})_t\|&\le&(2d\ep)^{p-q}\,\Ga^{2C}_t
\Big(\st\|\De_i^n X\|^{q}+\sum_{s\le t} \|\De X_s\|^{q} \Big)\\
&\le& 2d^{p-1}(2d\ep)^{p-q}\,\Ga^{2C}_t
\sum_{j=1}^d \sum_{s\le t} |\De X_s^j|^q.
\end{eqnarray*}
Since $\sum_{s\le t} |\De X_s^j|^q<\infty$, by letting $\ep\to0$
we conclude
$$\limsup_n \, \sup_{s\le t}~ |Z^n(f\Psi_\ep)_s|~=0,$$
which ends the proof in the case where $X$ is bounded.

The general case is deduced by a classical method of "localization", for
which we refer to Section 3 of \cite{BGJPN} for details.\qed
\vskip5mm

\noindent{\bf Proof of Theorem \ref{t2}:}
We use the previous notation, with $Z^n(f)$ is as in (\ref{ep37}) and
$D(f)$ as in (\ref{e15}). Recalling that (\ref{e16}) and (\ref{e15})
give the same process $D(f(1-\Psi_\ep))$, we still have (\ref{ep40}),
and it is thus enough to prove that:
\bee\label{e1}Z^n(f\Psi_{\ep})~\longrightarrow^{u.c.p.}~0. \eee

Set  $f_{\ep}:=f\Psi_{\ep}$.
By the hypotheses on $f$,  the function $f_{\ep}$  is
$C^2$ in $x$ if $\ep$ is small enough. We  then apply  lemma \ref{l2} to each
$f_{\ep}((i-1)\De_n,X_{(i-1)\De_n},\De_i^nX)$, which gives
$Z^n(f_{\ep})_t ~ = ~ \sum_{l=1}^3 Z^n(f_{\ep},l)_t$ where, with
the notation $Y^n_s=X_s-X_{(i-1)\De_n}$ and $\phi^n(s):=(i-1)\De_n$  for
$s\in((i-1)\De_n,i\De_n]$, we have
$$ \begin{array}{lll} Z^n(f_{\ep},1)_t & = & \sum_{j=1}^d
\int_0^{[t/\De_n]\De_n}
\left(\frac{\partial f_{\ep}}{\partial
x^j}(\phi^n(s),X_{\phi^n(s)},Y^n_s)-\frac{\partial f_{\ep}}{\partial
 x^j}(s-,X_{s-},0)\right)\,dX^j_s,\\[4mm]
 Z^n(f_{\ep},2)_t &  = &  \frac1{2}\sum_{j,k=1}^d\int_0^{[t/\De_n]\De_n}
\left(\frac{\partial^2 f_{\ep}}{\partial x^j \partial
x^k}(\phi^n(s),X_{\phi^n(s)},Y^n_s)\right.\\[3mm]
 & & \left.- \frac{\partial^2
f_{\ep}}{\partial
 x^j \partial x^k}(s-,Z_{s-},0)\right)\,d\langle X^{c,j},X^{c,k}
 \rangle_s,\\[4mm]
 Z^n(f,3)_t & = & \sum_{s\le t}
\Big(f_{\ep}(\phi^n(s),X_{\phi^n(s)},Y^n_s)-f_{\ep}(s-,X_{s-},\De X_s)\\[3mm]
  & & -
 f_{\ep}(\phi^n(s),X_{\phi^n(s)},Y^n_{s-}) -\sum_{j=1}^d \De X_s^j
 \Big(\frac{\partial f_{\ep}}{\partial x^j }(\phi^n(s)
 ,X_{\phi^n(s)},Y^n_{s-})\\[3mm]
  & & -\frac{\partial f_{\ep}}{\partial x^j
  }(s-,X_{s-},0)\Big)\Big).
\end{array} $$
Observe now that  $ \frac{\partial f_{\ep}}{\partial x^j}(\phi^n(s)
,Z_{\phi^n(s)},Y^n_{s-})\to \frac{\partial f_{\ep}}{\partial
x^j}(s-,Z_{s-},0)$.
Since $\frac{\partial f_{\ep}}{\partial x^j }$ is dominated by a locally
bounded processes, Lebesgue's theorem gives:
$$Z^n(f_{\ep},1)~\rightarrow^{u.c.p.}~0.$$
The proof of $Z^n(f_{\ep},j)\rightarrow^{u.c.p.}0$ for $j=2,3$ is
similar, and we thus have (\ref{e1}).

\subsection{Proof of Theorem \ref{t3}}

Let us start by strengthening  the hypothesis $(N_0)$ :

\vspace{.5cm}
\ni{\bf Hypothesis $(LN_0)$:} $(N_0)$ is satisfied,
and  the processes $b_s,\si_s,\int_{\R}(1\land
||\de(\om,s,y)||^2)\,F(dy)$ and  $X_s$ are bounded by a
constant.\qed
\vspace{.5cm}

We also  suppose that the process $\Ga$ which intervenes in
(\ref{er1}) is uniformly bounded. Below, we denote all constants
by $K$. Set
\bee\label{ep25} \be_i^n=\si_{(i-1)\De_n}\frac{\De_i^n W}{\sqrt{\De_n}}.
\eee

\begin{lem}\label{L44}
Suppose  $(LN_0)$ satisfied and  $f$ optional, satisfying $(K(\R))$  and
at most with polynomial growth. Then
\bee\label{e18}\De_n\st\,E\left\{f\left((i-1)\De_n,X_{(i-1)\De_n},
\be_i^n\right)|\f_{(i-1)\De_n}\right\}\,
\stackrel{u.c.p.}{\longrightarrow}\,\int_0^t H_{s-}\,ds,\eee
when $n\to\infty$, where
$ H_s=\int_0^t\rho_{\si_{s}}\left(f(s,X_{s},.)\right)\,ds$.
\end{lem}

\ni{\bf Proof:~} The left side of (\ref{e18}) is almost surely
equal to $\De_n\st H_{(i-1)\De_n}$.
This is a Riemann sum which therefore converges to $\int_0^t H_{s-}\,ds$
locally uniformly in $t$, because $H$ is a c\`adl\`ag process.\qed

\begin{lem}\label{L45}
Let  $f$ be optional, locally equicontinuous
in $x$ and    with  at most $p$-polynomial growth. Assume further that
$X$ satisfies $(LN_0)$ and either is continuous or $p<2$.
Then
$$\De_n\st E\left(\left\| f\Big((i-1)\De_n,X_{(i-1)\De_n},\frac{\De_i^n
X}{\sqrt{\De_n}}\Big)-f\Big((i-1)\De_n,X_{(i-1)\De_n},\be_i^n\Big)\right\|
\right)\rightarrow0.$$
\end{lem}
{\bf Proof:} We  reproduce the  proof of Lemma 4.4 (2) of [4]
with some relevant changes. For any $A,\,T,\,\ep >0$, we define
the variables
$$G_{T}(\ep,A)~=~
\sup_{s\le T;\,||x||\le A;
  ||z||\le K;\,||y||\le \ep}~\|f(s,z,x+y)-f(s,z,x)\|$$
$$ \ze_i^n ~= ~  \Big\| f\left((i-1)\De_n,X_{(i-1)\De_n},\frac{\De_i^n
X}{\sqrt{\De_n}}\right)-f((i-1)\De_n,X_{(i-1)\De_n},\be_i^n)\Big\|.$$
Then
\bee\label{ep2}\|\ze_i^n\|\,\le~
G_{t}(\ep,A)+\|\ze_i^n\|\left(1_{\left\{\|\be_i^n\|>A\right\}}+1_{\left\{
\left\|\De_i^n X/\sqrt{\De_n}-\be_i^n\right\|\,>\,\ep\right\}}\right).\eee

Let $q$ be a real such that $q>p$ if $X$ is continuous and $q=2$ if not.
Then (\ref{er1}) with $\Ga$ a constant yields for all $B>1$:
$$\|f(\om,s,z,x)\| \,\le \,
    K\,\phi(z)\Big(B^{p-q}\|x\|^q+B^p\Big). $$
Also  under $(LN_0)$ one knows that:
$$\E\left\{\left\|\De_i^n
X/\sqrt{\De_n}-\be_i^n\right\|^q+\|\be_i^n\|^q\right\}~\le~K.$$
Hence by (\ref{ep2}):
\begin{eqnarray*} \|\ze_i^n\| & \le &
G_{t}(\ep,A)+KB^p\left(1_{\{\|\be_i^n\|>A\}}+1_{\left\{\left\|\De
X_i^n/\sqrt{\De_n}-\be_i^n \right\|>\ep
\right\}}\right)\\
 & & +KB^{p-q}\left(\|\be_i^n\|^q+\left\|\De_i^n X/\sqrt{\De_n}-
 \be_i^n\right\|^q\right).\end{eqnarray*}
It follows that
\begin{eqnarray} \De_n\st\E\{\|\ze_i^n\|\} & \le &
t\left(\E\{G_t(\ep,A)\}~+~\frac{KB^p}{A}~+~KB^{p-q}\right)\nonumber\\
  & & \hskip-2cm+\,KB^p\ep^{-2}\De_n \st
\E\Big\{1\land\Big\|\De_i^n X/\sqrt{\De_n}-\be_i^n\Big\|^2
\Big\}.\label{ep3}\end{eqnarray}
Next by lemma 4-1 of (\cite{JJ})
$$\De_n \st \E\Big\{1\land \Big\|\De_i^n
X/\sqrt{\De_n}-\be_i^n\Big\|^2
\Big\}\,\longrightarrow\,0.$$ Then
coming back to (\ref{ep3}) and letting successively $n\to\infty,~
\ep\to 0,~A\to\infty$ and $B\to\infty$, we obtain the result. \qed
\vskip5mm

\noindent\bf Proof of Theorem \ref{t3}: \rm
We first prove  the theorem under the stronger assumptions
(LN$_0$) and $\Ga_t$ in (\ref{er1}) bounded. Set
$$U'^n_t:=\,\De_n\st f\left((i-1)\De_n,X_{(i-1)\De_n},\frac{\De_i^n
X}{\sqrt{\De_n}}\right)\,-\,\int_0^t\rho_{\si_{s}}\left(f(s-,X_{s},.)
\right)\,ds.$$
Then $U'^n_t=\sum_{j=1}^3 U'^n_t(j),$ where
\begin{eqnarray*} U'^n_t(1) & = & \De_n\st \Big(f\Big((i-1)
\De_n,X_{(i-1)\De_n},\frac{\De_i^n
X}{\sqrt{\De_n}}\Big) -\,f\left((i-1)\De_n,X_{(i-1)\De_n},\be_i^n\right)\Big),\\[2mm]
 U'^n_t(2)& = & \De_n \st
\Big(f\left((i-1)\De_n,X_{(i-1)\De_n},\be_i^n\right)\\[1mm]
 & & -\,
\E\left\{f\left((i-1)\De_n,X_{(i-1)\De_n},\be_i^n\right)|\f_{(i-1)\De_n}
\right\}\Big),\\[2mm]
U'^n_t(3) & = &  \De_n\st \E\left\{f\left((i-1)\De_n,X_{(i-1)\De_n},\be_i^n
\right)|\f_{(i-1)\De_n}\right\}  -\int_0^{t}\rho_{\si_{s}}\left(f\left(s-,X_s,.\right)\right)\,ds.
\end{eqnarray*}

Observe first that  $U'^n_t(2)$  is a martingale with respect  to the
filtration $(\f_{[t/\De_n]\De_n})_{t\ge 0}$, and its predictable quadratic
variation is given by:
\begin{eqnarray*} \langle U'^n(2)\rangle_t  &  = &   \De_n^2 \st
\Big(\E\Big\{f((i-1)\De_n,X_{(i-1)\De_n},\be_i^n)^2|\f_{(i-1)\De_n}\Big\}\\
 & & -\Big(\E\Big\{f((i-1)\De_n,X_{(i-1)\De_n},\be_i^n)\,|\f_{(i-1)\De_n}
 \Big\}\Big)^2\Big),\end{eqnarray*}
which satisfies  $\langle U'^n(2)\rangle_t\,\le\,Kt\De_n$.
It follows by  Doob's inequality, that  $U'^n_t(2)\stackrel{u.c.p.}
{\rightarrow}0$.
We have the same results for $U'^n_t(1)$ and $U'^n_t(3)$, respectively by
lemma \ref{L45} and \ref{L44}.

At this stage the theorem is proved under the stronger
assumptions announced at the beginning of the proof, and as said in
Theorem \ref{t1},
 the general case is obtained by a classical localization method.

\section{Proof of the central limit theorems}
\subsection{Proof of theorem \ref{t4}}

We  start again by strengthening our hypotheses:
\vspace{0.5cm}

\ni{\bf Hypothesis $(LN_1)$:}  $(N_1)$ is satisfied, and  the processes
$b$, $\si$ and  $X$ are bounded. The functions $\ga_k=\ga$ do not
depend on $k$ and are bounded.\qed
\vspace{0.5cm}

\ni{\bf Hypoth\`ese $(LM_1)$:} We have $(M_1)$ and the process
 $\Ga$  is bounded.\qed
\vspace{0.5cm}

Under $(LN_1)$, we have:  \bee\label{ep4}X_t = X_0  +
  \int_0^t b'_s\,ds+\int_0^t \si_s\,dW_s+\int_0^t\int_{\R}
  \de(s,y)(\umu-\unu)(ds,dy),\eee
where \bee\label{ep5} b'_s:= b_s + \int_{\R} h'(\de(s,y))\,F(dy).\eee

For $\ep>0$, set:
\bee\label{eoa1}E=\left\{\,y\in\R,~\ga(y)>\ep\,\right\}~~\mbox{and}~~
N_t=1_E\star\umu_t,\eee
and let  $T'_1,\cdots,T'_p,\cdots$ be the successive  jump times of $N$.

We state two important  lemmas, the first of which is due to Jacod and
Protter (Lemma 5.6 of \cite{JP}),
and the second  one is Lemma 5.9 of \cite{JJ}.

\begin{lem}\label{LO2}
Suppose $(LN_1)$ satisfied, and for each  $T'_p$, denote by $i^n_p$ the
integer such that
$(i^n_p-1)\De_n<T'_p\le i^n_p\De_n$. Then the sequence of random variables
$$\frac1{\sqrt{\De_n}}\left(\si_{(i^n_p-1)\De_n}(W_{T'_p}-
W_{(i^n_p-1)\De_n})~,~
 \si_{T'_p}(W_{i^n_p\De_n}-W_{T'_p})\right)_{p\ge 1}$$ converges stably in
 law to
  $$ \left(\sqrt{\kappa_p}\,\si_{T'_{p-}}U_p~,~
\sqrt{1-\kappa_p}\,\si_{T'_p}U'_p\right)_{p\ge1},$$
where $U_p$ is such that
$U^t_p=(U^1_p,\cdots,U^m_p)$ and $U'^t_p=(U'^1_p,\cdots,U'^m_p)$.
\end{lem}

\begin{lem}\label{LO1}
Under the assumptions of lemma \ref{LO2}, on has:
$$\frac1{\sqrt{\De_n}}\, \left((X_{i^n_p\De_n}-X_{T'_p}-
\si_{T'_p}(W_{i^n_p\De_n}-W_{T'_p})\right)\,\stackrel{\PP}{\longrightarrow}
\,0,$$
$$\frac1{\sqrt{\De_n}}\,\left(X_{T'_p-}-X_{(i^n_p-1)\De_n}-
\si_{(i^n_p-1)\De_n}(W_{T'_p}-W_{(i^n_p-1)\De_n})\right)\,
\stackrel{\PP}{\longrightarrow}\,0.$$
\end{lem}
\vspace{.4cm}

We are now ready to give the proof of the theorem.

The processes
\bee\label{ep6} W^n(f)=\frac1{\sqrt{\De_n}}\left(V^n(f)_t-
\sum_{s\le [t/\De_n]\De_n} f(s-,X_{s-},\De X_s)\right) \eee
satisfy $ W^n(f)=W^n(f,1)+W^n(f,2),$ where
\begin{eqnarray*} W^n(f,1)_t & = &
\frac1{\sqrt{\De_n}}\left(V^n(f)_t\,-
\sum_{s\le [t/\De_n]\De_n}f(\phi^n(s),X_{s-},\De X_s) \right)
\\[2mm]
W^n(f,2)_t & = & \frac1{\sqrt{\De_n}}\sum_{s\le [t/\De_n]\De_n}\left(
f(\phi^n(s),X_{s-},
\De X_s)-f(s-,X_{s-},\De X_s)\right),
\end{eqnarray*}
($\phi^n(s)$ is like in the previous section).
(\ref{e7}) yields $W^n(f,2)~\stackrel{u.c.p.}{\longrightarrow}~0$,
and for all $\ep>0$ we have
\bee\label{ep7} W^n(f,1)=W^n(f(1-\Psi_{\ep}),1)+W^n(f\Psi_{\ep},1),\eee
where $\Psi_\ep$ is as in  (\ref{e8}). Then the rest of the
proof of Theorem \ref{t4} is divided in three  steps.
\vspace{.5cm}

\ni{\bf Step 1:}
Here we study  the convergence of the process
$W^n(f(1-\Psi_{\ep}),1)$.
By subsection 3.1 of \cite{JJ}, for $n$ large enough one has:
\begin{eqnarray*} W^n(f(1-\Psi_{\ep})) & = &
\frac1{\sqrt{\De_n}} \sum_{p:~T'_p\le [t/\De_n]\De_n}\left(
f(1-\Psi_{\ep})((i^n_{p}-1)\De_n,X_{(i^n_{p}-1)\De_n},\De_i^n X)\right.
\\& & \left. -~f(1-\Psi_{\ep})((i_p^n-1)\De_n,X_{T'_p-},\De X_{T'_p})
\right),\\
 & =  & \frac1{\sqrt{\De_n}} \sum_{p:~T'_p\le [t/\De_n]\De_n}
\left(\sum_{j=1}^d\left( \De_{i_p^n}^n X^j-\De X^j_{T'_p} \right)\right.\\
 & & \ti\frac{\partial f(1-\Psi_{\ep})}{\partial x_j}((i^n_{p}-1)\De_n,
\overline{X}'^{n}_p,\overline{X}^n_p) \\
& & \left.+ \sum_{j=1}^d\left( X^j_{(i_p^n-1)\De_n}-X^j_{T'_p-}\right)
\frac{\partial f(1-\Psi_{\ep})}{\partial z_j}
((i^n_{p}-1)\De_n,\overline{X}^{'n}_p,\overline{X}^n_p) \right),
  \end{eqnarray*}
where $(\overline{X}'^{n}_p,\overline{X}^n_p)$ is between
$(X_{(i^n_{p}-1)\De_n},\De_i^n X)$ and
$(X_{T'_p-},\De X_{T'_p})$.
Then by lemma \ref{LO1} and \ref{LO2}, $W^n(f(1-\Psi_{\ep})$  converge
stably in law to the process
\begin{eqnarray*} F'(f(1-\Psi_{\ep}))_t &:=  & \sum_{p:~T'_p\le t}
\sum_{j=1}^d \sum_{k=1}^m \left( \left(
\sqrt{\ka_p} \,\si^{j,k}_{T'_p-}U_p^k+\sqrt{1-\ka_p} \,\si^{j,k}_{T'_p}
{U}_p^{'k}\right) \right.\\
  & & \ti \frac{\partial f(1-\Psi_{\ep})}{\partial x_j}
  (T'_p-,X_{T'_p-},\De X_{T'_p})\\
 & & \left.  - \sqrt{\ka_p}\,\si^{j,k}_{T'_p-}U_p^k\frac{\partial
 f(1-\Psi_{\ep})}{\partial z_j}(T'_p-,X_{T'_p-},\De X_{T'_p})\right),
 \end{eqnarray*}
which has the same $\f$-conditional law than the process
$F(f(1-\Psi_\ep))$ associated with the function $f(1-\Psi_\ep)$ by
(\ref{ep15}).
\vspace{.5cm}

\ni{\bf Step 2:} Here we show that
\bee\label{ep10} F(f(1-\Psi_{\ep})) ~\stackrel{u.c.p.}{\longrightarrow}
~F(f) \quad \mbox{as}~\ep\to~0. \eee
Recall the process $C(f)$ defined in (\ref{e21}),
and set $f\Psi_{\ep}=f_{\ep}$.
Under $(LN_1)$ there exists a process $A$ such that:
$$\forall\,T>0,~~ C(f_{\ep})_T\le A_T,\quad \mbox{and}\quad
\E(A_t)<\infty.$$
Since $C(f_{\ep})_T \to 0$  when $\ep\to 0$,  by Lebesgue's convergence
theorem we have $\E(C(f_{\ep})_T)\to0$.
Furthermore by lemma \ref{l3}, the process  $F(f_{\ep})_t$ is a locally
square integrable martingale
and Doob's inequality yields that:
$$\widetilde{\PP}\left( \sup_{t\le T}\, \| F(f_{\ep})\|\, >\,\eta \right) \le
\frac4{\eta}~\widetilde{\E}\left(F(f_{\ep})^2_t\right)=
\frac4{\eta}\E\left(C(f_{\ep})_T\right)$$
hence $F(f\Psi_{\ep})\stackrel{u.c.p.}{\longrightarrow} 0$ when
$\ep\to0$. Since $F(f)= F(f(1-\Psi_{\ep}))+ F(f\Psi_{\ep})$, this implies
(\ref{ep10}).
\vspace{.5cm}

\ni{\bf Step 3:}
In this last step we show  that
\bee\label{ep11} \lim _{\ep\to 0}~\limsup_n ~\PP\left\{ \sup_{t\le T}
\|W^n(f\Psi_{\ep},1)\|\,>\,\eta\right\} ~=~0,\quad \forall\,\eta,\,T\,>0.
\eee

Using  It\^o's formula of lemma \ref{l2}, in a similar way than  in the
proof of theorem \ref{t2}, we have
$W^n(f_{\ep},1)=\sum_{l=1}^5 W^n(f_{\ep},1,l)$, where
\begin{eqnarray*}
W^n(f_{\ep},1,1) & = & \frac{1}{\sqrt{\De_n}}~
\sum_{j=1}^d\int_{0}^{[t/\De_n]\De_n} {b'}^j_s\frac{\partial f_{\ep}}
{\partial x_j}\left(\phi^n(s),X_{\phi^n(s)},Y^n_{s-}\right)
\,ds, \\
W^n(f_{\ep},1,2)_t & = &
\frac1{2\sqrt{\De_n}}\sum_{j,j'=1}^d\sum_{k=1}^m
\int_0^{[t/\De_n]\De_n}\si^{j,k}_{s}\si^{j',k}_{s}
 \frac{\partial^2f_{\ep}}{\partial x_j\partial x_{j'}}
\left(\phi^n(s),X_{\phi^n(s)},Y^n_{s-}\right)\,ds,\\
W^n(f_{\ep},2,3)_t & = & \frac1{\sqrt{\De_n}}
\sum_{j=1}^d\sum_{k=1}^m\int_0^{[t/\De_n]\De_n}\si^{j,k}_s\frac{\partial f_{\ep}}{\partial x_j}
\left(\phi^n(s),X_{\phi^n(s)},Y^n_{s-}\right)
\,dW^k_s,\\
W^n(f_{\ep},1,4)_t & = & \frac1{\sqrt{\De_n}}\sum_{j=1}^d
\int_0^{[t/\De_n]\De_n}\int_0^{\R}\de^j(s,y)\frac{\partial f_{\ep}}{\partial x^j}
\left(\phi^n(s),X_{\phi^n(s)},Y^n_{s-}\right)
\,(\umu-\unu)(ds,dy),\\
W^n(f_{\ep},1,5)_t & = & \frac1{\sqrt{\De_n}}
\int_0^{[t/\De_n]\De_n}\int_0^{\R}\Big( f_{\ep}
\left(\phi^n(s),X_{\phi^n(s)},Y_{s-}^n+
\de(s,y)\right) \\
& & \hskip-2cm -f_{\ep}\left(\phi^n(s),X_{\phi^n(s)},Y^n_{s-}\right)
-\sum_{j=1}^d\de^j(s,y)\frac{\partial f_{\ep}}{\partial x_j}
\left(\phi^n(s),X_{\phi^n(s)},Y^n_{s-}\right)\\
& &  \hskip2cm- f_{\ep}(\phi^n(s),X_{s-},\de(s,y))\Big)~\umu(ds,dy),
\end{eqnarray*}
with $Y^n_s$ and $\phi^n(s)$ as before.

Under $(LN_1)$ and $(LM_1)$,  we have:
\bee\label{ep14}\left. \begin{array}{lll} \E(||X_t-X_s||^p)  & \le &
K|t-s|^{p/2},~~\forall p\in[0,2].\\[2mm]
\sum_{j=1}^d\left\|\frac{\partial f_{\ep}}{\partial z_j}(s,z,x)\right\| &
\le  & \al_{\ep}
(||x||\land(2d\ep))^2,\\[2mm]
\sum_{j=1}^d\sum_{j'=1}^d\Big\|\frac{\partial^{2}
f_{\ep}}{\partial x_j\partial z_{j'}}(s,z,x)\Big\| & \le  & \al_{\ep}
(||x||\land(2d\ep)),\\[2mm]
\sum_{j=1}^d \sum_{j'=1}^d\Big\|\frac{\partial^{k_1+k_2} f_{\ep}}
{\partial x^{k_1}_j\partial x^{k_2}_{j'}}(s,z,x)\Big\| & \le  &
\al_{\ep}(||x||\land(2d\ep))^{3-(k_1+k_2)},
\end{array} \right\}
\eee
where $\al_{\ep}\to0$  when $\ep\to 0$, $k_1+k_2\in\{0,1,2\}$,
and $\frac{\partial^0(f_{\ep})}{\partial x_j^0}=f_{\ep}$. We also
have
$$\frac{\partial f_{\ep}}{\partial x_j}
\left(\phi^n(s),X_{\phi^n(s)},Y^n_{s-}\right)=
\sum_{j'=1}^d Y^{n,j'}_{s-}\frac{\partial^2 f_{\ep}}{\partial x_j
\partial x_{j'}}
\left(\phi^n(s),X_{\phi^n(s)}, \overline{Y}^n_s\right),$$
where $\overline{Y}^n_s$ belongs to the segment
joining $Y^n_{s-}$ and $0$, thus
\begin{eqnarray} &&\int_0^{[t/\De_n]\De_n}\E\left\{\left\|
\frac{\widehat{b}^{j}_s~Y^{n,j'}_{s-}}{\sqrt{\De_n}}
\frac{\partial^2f_{\ep}}{\partial x_j\partial x_{j'}}
\left(X_{\phi^n(s)},\overline{Y}^n_s\right)\right\|\right\}\,ds\nonumber\\
&&  \le~
K\int_0^{t}\left[\left(\E\left\{\frac{\|Y_{s-}^{n,j'}\|^2}{\De_n}\right\}
\right)^{1/2}
\left(\E\left\{\left\|\frac{\partial^2f_{\ep}}{\partial x_j\partial x_{j'}}
\left(\phi^n(s),Z_{\phi^n(s)},\overline{Y}^n_s\right)\right\|^2\right\}\right)
^{1/2}\right]\,ds\nonumber\\
& &\le~    K \int_0^{t}\left( \E\left\{\left\|\frac{\partial^2f_{\ep}}
{\partial x_j\partial x_{j'}}
\left(\phi^n(s),Z_{\phi^n(s)},\overline{Y}^n_s\right)\right\|^2\right\}
\right)^{1/2}\,ds.
\label{c3e34} \end{eqnarray}
Since $\frac{\partial^2f_{\ep}}{\partial x_j\partial x_{j'}}(\om,s,z,0)
=0$, and $\frac{\partial^2f_{\ep}}{\partial x_j\partial x_{j'}}(\om,s,z,x)$
satisfies $(K(V))$, one deduces by Lebesgue's theorem that
(\ref{c3e34}) converge to $0$,  and thus
$$W^n(f_{\ep},2,1)~\stackrel{u.c.p.}{\longrightarrow}~0.$$
Similarly we show that $$W^n(f_{\ep},2,2)~
\stackrel{u.c.p.}{\longrightarrow}~0.$$

Next, the processes
$$\frac1{\sqrt{\De_n}}\int_0^{[t/\De_n]\De_n}\frac{\partial f_\ep}{\partial
x_j}\left(\phi^n(s),X_{\phi^n(s)},Y^n_{s-}\right)\si^{j,k}_sdW^k_s$$
are martingale with respect to the filtration $(\f_{[t/\De_n]\De_n})$, hence
 by Doob's inequality  and (\ref{ep14}) one has:
$$\PP\left\{
\sup_{t\le T}~\left\|\frac1{\sqrt{\De_n}}\int_0^{[t/\De_n]\De_n}
\frac{\partial f_{\ep}}{\partial x_j}
\left(\phi^n(s),X_{\phi^n(s)},Y^n_{s-}\right)\si^{j,k}_sdW^k_s\right\|~
>\eta\right\}$$
$$\le ~ \frac{1}{\eta^2\De_n}\int_0^{T}\E\left\{
\left\|\frac{\partial f_{\ep}}{\partial
x}\left(\phi^n(s),X_{\phi^n(s)},Y^n_{s-}\right)\si^{j,k}_s\right\|^2\right\}
\,ds  ~ \le ~ \frac{KT\al_{\ep}^2}{\eta^2},$$
and
$$\lim_{\ep\to 0}~\limsup_n~\PP\left\{\sup_{t\le
T}\left\|W^n(f_{\ep},1,3)_t\right\|~>\eta\right\}~=~0.$$
Similarly, we  have:
$$\lim_{\ep\to 0}~\limsup_n~\PP\left\{ \sup_{t\le
T}\left\|W^n(f_{\ep},1,4)_t\right\|~>\eta\right\}~=~0.$$

Now under  $(LM_1)$, separating the cases where  $||x||\le ||x'||$
and  $||x'||\le ||x||$, one  shows that:
$$\Big\| f_{\ep}(\om,s,z_1,x+x')-f_{\ep}(s,z_1,x')-
\sum_{j=1}^d x_j\frac{\partial f_{\ep}}{\partial x_j}(\om,s,z_1,x') -
f_{\ep}(\om,s,z_2,x)\Big\|$$
$$\le~K\al_{\ep}||x||^2\left(~||z_1-z_2||~+~||x'||~\right).$$
Then $\PP\left( \sup_{t\le T} |W^n(f_\ep,1,5)_t| >\eta\right)$
is smaller than
\begin{eqnarray*}
&&  ~\frac1{\eta}~\E\Big\{\int_0^{[t/\De_n]\De_n}\int_{\R}~\|~
f_{\ep}(\phi^n(s),X_{\phi^n(s)},Y^n_{s-}
+\de(s,y))\\[2mm]
& &-~ f_{\ep}(\phi^n(s),X_{\phi^n(s)},\de(s,y))
-~f_{\ep}(\phi^n(s),X_{\phi^n(s)},Y^n_{s-})\\
 & &  -~ \de(s,y)\frac{\partial f_{\ep}}{\partial x}(\phi^n(s),
 X_{\phi^n(s)},Y^n_{s-})~\| \umu(ds,dy)~\Big\}\\
&&\qquad \le  K\al_{\ep}\left(\int_0^t\E\left\{\frac{||Y_{s-}^n||+
||Z_{\phi^n(s)}-Z_{s-}||}{\sqrt{\De_n}}\right\}~ds\right)~\le~
K\al_{\ep}t \end{eqnarray*}
and thus
$$\lim_{\ep\to 0}~\limsup_n~\PP\left\{
\sup_{t\le
T}\left\|W^n(f_{\ep},1,5)_t\right\|~>\eta\right\}~=~0.$$
This ends the  proof under the reinforced assumptions $(LN_1)$ and $(LM_1)$.
One finishes
the proof by a classical localization procedure.

\subsection{Proof of theorem \ref{t5} and \ref{t6}.}\label{s2}

As for the previous proofs, we first strengthen the hypotheses, and
thanks to Remark \ref{r1}, we adopt the form (\ref{e9}) for $\si$.
\vspace{.5cm}

\ni{\bf Hypothesis $(LN_2(s))$:}
We have $(N_2(s))$ and the processes $b_s$, $\widetilde{b}_s$,
$\widetilde{\si}_s$, $\widetilde{v}_s$,
$\De\si_s$, $\int_{\R} (1\land ||\widetilde{\de}(s,y)||^2)\,F(dy)$
are bounded. The functions $\ga_k=\ga$ do not depend on $k$ and are also
bounded.  \qed
\vspace{.5cm}

We denote $(LM_2)$ (resp. $(LM'_2)$)
the hypothesis $(M_2)$ (resp. $(M'_2)$) with the
additional condition that the process $\Ga$ is bounded.

Under $(LN_2(s))$ with  $s\le 1$, $X$ can be write as:
\bee\label{ep19} X_t~=~X_0+\int_0^t b'_s\,ds+\int_0^t\si_s\,dW_s+
\int_{\R}\int_0^t \de(s,y)\,\umu(ds,dy),\eee
where $b'_s=b_s-\int_{\R} h(\de(s,y))\,F(dy)$,
and  under  $(L_2(2))$ the process  $\si$ is writen:
\bee\label{ep20}\si_t ~ = ~ \si_0+\int_0^t\widetilde{b}'_s\,ds
+\int_0^t\widetilde{\si}_s\,dW_s+
\int_0^t \widetilde{v}\,dV_s+ \int_0^t\int_{\R}\widetilde{\de}(s,y)\,(ds,dy),
\eee
with $\widetilde{b}'_s=\widetilde{b}_s+\int_{\R} k
'\left(\widetilde{\de}(s,y)\right)F(dy).$

Let us now give some useful lemmas.

\begin{lem}\label{l8}
Suppose $(LN_2(2))$ satisfied and assume that  $f$ is optional, locally
equicontinuous in $x$ and  at most with  $p$-polynomial growth. If further,
either  $X$ is  continuous or $p<1$, then:
\bee\label{ep21}
\De_n\st\,E\left(\left\| f\left((i-1)\De_n,X_{(i-1)\De_n},\frac{\De_i^n
X}{\sqrt{\De_n}}\right)-f((i-1)\De_n,X_{(i-1)\De_n},\be_i^n)\right\|^2
\right)\to0.
\eee
\end{lem}

\ni{\bf Proof:} The proof of this lemma is the same as for
Lemma \ref{L45},  the condition $p<1$ come in because of the
the square in (\ref{ep21}). \qed
\vskip5mm

Set
\begin{eqnarray}\label{ep43}
U^n_t & := & \sqrt{\De_n}\st\E\Big\{\Big(f\Big((i-1)\De_n,X_{(i-1)\De_n},
\frac{\De_i^n
X}{\sqrt{\De_n}}\Big) \nonumber\\
& &
-f\left((i-1)\De_n,X_{(i-1)\De_n},\be_i^n\right)\Big)\,|\f_{(i-1)\De_n}\Big\}
\end{eqnarray}

\begin{lem}\label{l4}
Suppose  $(LN_2(2))$ and $(M_2)$ satisfied and  $X$ continuous.
Assume further that  one of the following two conditions is satisfied:
\begin{itemize}
\item[{\bf A.}]The application  $x\mapsto f(\om,s,z,x)$ is even in $x$.
\item[{\bf B.}] We have $b'=0$ and $\widetilde{\si}=0$.
\end{itemize}
Then $U^n\,\stackrel{u.c.p.}{\longrightarrow}\,0$.
\end{lem}

\ni{\bf Proof:} {\bf A)} Set  \bee\label{c4e35}
L_i^n:=f\left((i-1)\De_n,X_{(i-1)\De_n},\frac{\De_i^n
X}{\sqrt{\De_n}}\right)-f\left((i-1)\De_n,X_{(i-1)\De_n},\be_i^n\right).
\eee Then $L_i^n=L_i^{'n}+ L_i^{''n}$, where
\begin{eqnarray*}
L'^n_i &  = &  \sum_{j=1}^d\left(\frac{\partial f}{\partial
x_j}\left((i-1)\De_n,X_{(i-1)\De_n},\bar{\ga}^n_i\right)~-~\frac{\partial f}
{\partial x_j}\left((i-1)\De_n,X_{(i-1)\De_n},\be_i^n\right)\right)\nonumber\\
& & \ti \left(\frac{\De_i^n X^j}{\sqrt{\De_n}}-\be_i^{n,j} \right),\\
L''^n_i &  = &  \sum_{j=1}^d\left(\frac{\De_i^n X^j}
{\sqrt{\De_n}}-\be_i^{n,j}\right)\frac{\partial f}{\partial x_j}
\left((i-1)\De_n,X_{(i-1)\De_n},\be_i^n\right),\end{eqnarray*}
for some  random variable $\bar{\ga}_i^n$ between
$\frac{\De_i^n X}{\sqrt{\De_n}}$ and $\be_i^n$.
 For any $\ep,~A~>0$,  set
\begin{eqnarray*} G_{t}^A(\om,\ep) & =  & \sup_{
s\le t;\,\|y\|\le \ep;\,z\in{\cal K};\,\|x\|\le A}
\left\{~\sum_{j=1}^d\,\Big\|\frac{\partial f}{\partial x}
(\om,s,z,x+y)-\frac{\partial
f}{\partial x}(\om,s,z,x)\,\Big\|\right\}.
\end{eqnarray*}
Then:
\begin{eqnarray*}
\|L'^n_i\| & \le &  K\left(G_t^A(\ep)\,+\,\left( 1+\|\be_i^n\|^p+
\left\|\frac{\De_i^n X}{\sqrt{\De_n}}-\be_i^n\right\|^p\right)
\left(\frac{||\be_i^n||}{A}+ \frac{\left\|\frac{\De_i^n X}
{\sqrt{\De_n}}-\be_i^n\right\|}{\ep}\right)\right)
\\ & & \ti\left\|\frac{\De_i^n
X}{\sqrt{\De_n}}-\be_i^n\right\|.
\end{eqnarray*}

Next, under the assumption $(N_2(s))$ (in particular the properties of
$\si$), one shows that for all $q\geq2$:
\bee\label{c4e29}
\E(||\be_i^n||^q)\le K, ~~E\left\{\left\|\De_i^n X/
\sqrt{\De_n}-\be_i^n\right\|^q\right\}\le K \De_n.
\eee
Thus  by a repeated  use of H\"older inequality:
\bee\label{c4e30}\sqrt{\De_n}\st E\{\|L'^n_i\|\}~\le~
Kt\left[\left(E\left\{\left(G_t^A(\ep)\right)^2\right\}\right)^{1/2}+
\frac{\De_n^{1/4}}{\ep}~+~\frac1{A}\right].\eee
Letting successively $n\to\infty$, then $\ep\to0$ and then $A\to\infty$,
we obtain
\bee\label{c4e31} \sqrt{\De_n}\st E\{\|L'^n_i\|\}\longrightarrow0.\eee

 Let us now turn to $L''^{n}_i$. Under  $(LN_2(s))$ we have: $\frac{\De_i^n
X}{\sqrt{\De_n}}-\be_i^n=\widetilde{\xi}_i^n+\widehat{\xi}_i^n,$ where
\begin{eqnarray*}
\widehat{\xi}_i^n & = & \frac1{\sqrt{\De_n}}
\left(\int_{(i-1)\De_n}^{i\De_n}(b'_s-b'_{(i-1)\De_n})\,ds+
\int^{i\De_n}_{(i-1)\De_n}\Big(\int^s_{(i-1)\De_n} \widetilde{b'}_u\,du\right.
\\
 & & + \left.\int^s_{(i-1)\De_n}(\widetilde{\si}_u-
\widetilde{\si}_{(i-1)\De_n})\,dW_u\Big)\,dW_s\right),\\
\widetilde{\xi}^n_i & = & \sqrt{\De_n}b'_{(i-1)\De_n}
+\frac1{\sqrt{\De_n}}\int^{i\De_n}_{(i-1)\De_n}
\Big( ~\widetilde{\si}_{(i-1)\De_n}(W_s-W_{(i-1)\De_n})\\
 & &  +\int^{s}_{(i-1)\De_n}\int_{\R}\widetilde{\de}(u,y)(\umu-\unu)(du,dy)+
 \int^s_{(i-1)\De_n}\widetilde{v}_u\,dV_u
\Big)\,dW_s.\end{eqnarray*}

\ni{\bf 1) } Here we show that for any $j\in\{1,\cdots,d\},$
\bee\label{c4e22b} \E\left\{\widetilde{\xi}_i^{n,j}\frac{\partial f}{\partial
x_j}\left((i-1)\De_n,X_{(i-1)\De_n},\be_i^n\right)~|\f_{(i-1)\De_n}
\right\}~=~0.\eee
Since the function $x\to\frac{\partial f}{\partial x_j}(\om,s,z,x)$ is odd,
one  clearly has:
\bee\label{ep17}\E\left\{{b'}^j_{(i-1)\De_n}\frac{\partial f}{\partial
x_j}\left((i-1)\De_n,X_{(i-1)\De_n},\be_i^n\right)~|\f_{(i-1)\De_n}
\right\}\,=\,0,\eee
and for any $k,~k'\in\{1,\cdots,m\}$:
$$\E\left\{ \widetilde{\si}^{j,k,k'}_{(i-1)\De_n}
\Big(\int^{i\De_n}_{(i-1)\De_n}(W^{k'}_s-W^{k'}_{(i-1)\De_n})\,dW^{k}_s\Big)
\ti \right.$$
\bee\label{ep18}\left.\ti\frac{\partial f}{\partial x_j}
\left((i-1)\De_n,X_{(i-1)\De_n},\be_i^n\right)\,|\f_{(i-1)\De_n}
\right\}\,=\,0.\eee

Next consider  the $\si$-field:
$$\f'_{(i-1)\De_n}~=~\f_{(i-1)\De_n}\bigvee
\si(W_s-W_{(i-1)\De_n}:(i-1)\De_n\le s\le i\De_n).$$
\ni Since $W$  is independent of  $\umu$ and of $V$, for any $j,\,k$ as above
 one has:
$$ \E\left\{  \left( \int^{i\De_n}_{(i-1)\De_n} \left(
\int^{s}_{(i-1)\De_n}\int_{\R}\widetilde{\de}^{j,k}(u,y)(\umu-\unu)(du,dy)
\right)\,dW^k_s\right)\ti \right.$$
\bee\label{ep41}  \left. \ti \frac{\partial f}{\partial
x_j}\left((i-1)\De_n,X_{(i-1)\De_n},\be_i^n\right)~|\f_{(i-1)\De_n}
\right\}~=0,\eee
and for any $j'\in\{1,\cdots,l\}$: $$ \E\left\{
\Big( \int^{i\De_n}_{(i-1)\De_n} \Big(
\int^s_{(i-1)\De_n}\widetilde{v}^{j,k,j'}_u\,dV^{j'}_u
\Big)\,dW^k_s\Big)\ti\right.$$ \bee\label{ep26}\left. \ti\frac{\partial f}
{\partial
x_j}((i-1)\De_n,X_{(i-1)\De_n},\be_i^n)~|\f_{(i-1)\De_n}\right\}~=0.\eee

From   (\ref{ep17}), (\ref{ep18}), (\ref{ep41}) and (\ref{ep26}) we deduce
(\ref{c4e22b}).

\vspace{.5cm}

\ni{\bf 2)} In this step, we show  that for all $j\in\{1,\cdots,d\}$,
\bee\label{c4e22}
\sqrt{\De_n}\st\E\left\{\left\|\widehat{\xi}_i^{n,j}
\frac{\partial f}{\partial
x_j}\left((i-1)\De_n,X_{(i-1)\De_n},\be_i^n\right)\right\|~
|\f_{(i-1)\De_n} \right\}\longrightarrow~0.\eee
By   H\"older and  Doob inequalities  we have:
$$ \E\{||\,\widehat{\xi}_i^n\,||^2\}\,\le \,K\left(\De_n^3+
 \int_{(i-1)\De_n}^{i\De_n}\left(||b'_s-b'_{(i-1)\De_n}||^2
 +||\widetilde{\si}_s-
 \widetilde{ \si}_{(i-1)\De_n}||^2 \right)\,ds
\right).$$
 Since $\E\left\{\Big\|\frac{\partial f}{\partial
x}\left((i-1)\De_n,X_{(i-1)\De_n},\be_i^n\right)\Big\|^2\right\}~\le~K,$
it follows from a repeated use of Hold\"er inequality that
$$ \sqrt{\De_n}\st\E\left\{\left\|\widehat{\xi}_i^{n,j}
\frac{\partial f}{\partial
x_j}\left((i-1)\De_n,X_{(i-1)\De_n},\be_i^n\right)\right\|~|\f_{(i-1)\De_n}
\right\}~\le~Kt\De_n~+$$
$$+ Kt^{1/2}\left(
\E\left\{\int_{0}^{[t/\De_n]\De_n}\left(||b'_s-b'_{[s/\De_n]\De_n}||^2
+||\widetilde{\si}_s-\widetilde{\si}_{[s/\De_n]\De_n}||^2
\right)\,ds\right\}\right)^{1/2}.$$

Since $b'$ and $\widetilde{\si}$ have some continuity properties in $s$,
we deduce by Lebesgue theorem that the last quantity tends to $0$
when $n\to\infty$, hence (\ref{c4e22}).
\vskip2mm

\ni{\bf B)} The proof is the same than for
(A), except for the fact that we have (\ref{ep17}) and (\ref{ep18})
because $b'=\widetilde{\si}=0$.\qed
\vspace{.4cm}

We give now another version of  Lemma \ref{l4}, in  the
case where $X$ is discontinuous:

\begin{lem}\label{l7}
Suppose  $X$ satisfies $(LN_2(s))$  with $s\le 1$ and $f$ satisfies $(LM_2')$.
Assume further that either  $f(\om,s,z,x)$ is even in $x$ or ${b'}=\widetilde{\si}=0$.
Then $$ U^n\,\stackrel{u.c.p.}{\longrightarrow}\,0, \quad \mbox{when}~n\to\infty.$$
\end{lem}
\ni{\bf Proof:} Recall that under $(LN_2(s))$ with $s\le 1$,  $X$ is written
as in
 (\ref{ep19}).  Set
$$X'_t~:=~X_0~+~\int_0^t {b'}_s\,ds~+~\int_0^t \si_{s}\,dW_s.$$
Let  $(\ep_n)$ be a sequence such that:
$ \ep_n\in]0,1]$ and $\ep_n\to 0$ when $n\to\infty$,
and  set $E_n\,=\left\{x\in\R,\,\ga(x)> \ep_n \right\}.$
Then
\begin{eqnarray*}\frac{\De_i^n X }{\sqrt{\De_n}} & = &
\frac{\De_i^n X'}{\sqrt{\De_n}}~+~
\frac1{\sqrt{\De_n}}\int_{(i-1)\De_n}^{i\De_n}
\int_{E_n^c}\de(s,x)\,\umu(ds,dx).\\
 & & + \frac1{\sqrt{\De_n}}\int_{(i-1)\De_n}^{i\De_n}\int_{E_n}\de(s,x)\,
 \umu(ds,dx).\end{eqnarray*}

Set
\begin{eqnarray*} \ze_i^n (1) & := & \frac1{\sqrt{\De_n}}
\int_{(i-1)\De_n}^{i\De_n}\int_{E_n}\de(s,y)\,\umu(ds,dy), \\
\ze_i^n(2) & := & \frac1{\sqrt{\De_n}}
\int_{(i-1)\De_n}^{i\De_n}\int_{E^c_n}\de(s,y)\,\umu(ds,dy).\end{eqnarray*}
Then using  the notation (\ref{c4e35}), one has  $L_i^n=\sum_{j=1}^3 L_i^n(j)$, where

\begin{eqnarray*}
L_i^n(1)  =   f\Big((i-1)\De_n,X_{(i-1)\De_n},\frac{\De_i^n X}
{\sqrt{\De_n}}\Big)
-f\Big((i-1)\De_n,X_{(i-1)\De_n},\frac{\De_i^n X}{\sqrt{\De_n}}-\ze_i^n(1)
\Big),\\
L_i^n(2)  =   f\Big((i-1)\De_n,X_{(i-1)\De_n},\frac{\De_i^n X}{\sqrt{\De_n}}
-\ze_i^n(1)
\Big)-f\Big((i-1)\De_n,X_{(i-1)\De_n},\frac{\De_i^n X'}{\sqrt{\De_n}}\Big),\\
L_i^n(3) \, = \,f\Big((i-1)\De_n,X_{(i-1)\De_n},\frac{\De_i^n X'}
{\sqrt{\De_n}}\Big)-
f\Big((i-1)\De_n,X_{(i-1)\De_n},\be_i^n\Big).\end{eqnarray*}

The hypothesis $(LM'_2)$ gives the existence of a sequence of reals
$(K^m)$ such that
$$\|z\|\le m\,\Rightarrow\,\|f(\om,s,z,x_1)-
f(\om,s,z,x_1+x_2)\|~\le~K^m(1\land\|x_2\|).$$
Hence
$$ \sqrt{\De_n}\st\E\Big\{\,\|L_i^n(1)\|\,|\f_{(i-1)\De_n}\,\Big\}\,\le\,
K\sqrt{\De_n}\st\E\Big\{\,(1\land\|\ze_i^n(1)\|)\,|\f_{(i-1)\De_n}\,\Big\}.$$
By the inequality (5.9) of lemma 5.3 of \cite{JJ}, we deduce:
\bee\label{c4e26}\sqrt{\De_n}\st\E\left\{\left\|L_i^n(1)\right\|
\right\}~\le~Kt\De_n^{1/2}\ep^{-1}_n.\eee

Next, set $\te(y)= \int_{\{|\ga(x)|\le y\}}|\ga(x)|\,F(dx)$,
which goes to $0$ as $y\to0$. One has
\bee\label{c4e40}
\sqrt{\De_n}\st\E\{\|L_i^n(2)\|\}~\le~
K\sqrt{\De_n}\st\E\{\|\ze_i^n(2)\|\}~\le~Kt\te(\ep_n).\eee
Finally, lemma \ref{l7}  implies:
\bee\label{ep44} \sqrt{\De_n}\st\left\|\E\{L_i^n(3)~|\f_{(i-1)\De_n}\}
\right\|~\longrightarrow^{u.c.p.}~0,
\quad\mbox{when}~n\to\infty.\eee

By  (\ref{c4e26}),  (\ref{c4e40})  and (\ref{ep44}) we have:
\begin{eqnarray*} \sqrt{\De_n}\st\left\|\E\left\{L_i^n\,|\f_{(i-1)\De_n}
\right\}\right\| &  \le &
Kt\left(\De_n^{1/2}\ep^{-1}_n+\te(\ep_n)\right)\\  & & + \sqrt{\De_n}\st
\left\|\E\left\{L_i^n(3)~|\f_{(i-1)\De_n}\right\}\right\|. \end{eqnarray*}

Choosing   $\ep_n=(1\land\De^{1/4}_n)$,  we conclude:
$$ \sqrt{\De_n}\st\left\|\E\{L_i^n~|\f_{(i-1)\De_n}\}\right\|~
\longrightarrow^{u.c.p.}~0, $$
and this ends the proof. \qed
\vspace{.5cm}

Set now
\bee
U'^{n}_t
 =  \frac1{\sqrt{\De_n}}
\Big(\,\De_n\st\E\{f((i-1)\De_n,X_{(i-1)\De_n},\be_i^n)|\f_{(i-1)\De_n}\}
  -\int_0^t\rho_{\si_{s}}(f(s,X_{s},.))ds\Big). \label{c4e36}
\eee

\begin{lem}\label{L446}
If $X$ satisfies $(LN_2(2))$ and  $f$ satisfies $(LM_2)$, we have
$ U'^n\stackrel{u.c.p.}{\longrightarrow}0$.
\end{lem}

{\bf Proof:} We can assume without loss of generality that $f$ is
$1$-dimensional. We also write the proof when the dimensions of $X$ and
$\si$ are $1$, since the multidimensional case is more cumbersome
but similar to prove. We have
$U'^n_t=U'^n_t(1)+U'^n_t(2)+U'^n_t(3)$, where
\begin{eqnarray*}
U'^{n}_t(1) & = & \frac1{\sqrt{\De_n}}\st\int_{(i-1)\De_n}^{i\De_n}\Big(
\rho_{\si_{(i-1)\De_n}}(f((i-1)\De_n,X_{(i-1)\De_n},.))  \\
& &  -~ \rho_{\si_{s}}(f((i-1)\De_n,X_{s},.))\Big)\,ds,\\
U'^n_t(2) & = &   \frac1{\sqrt{\De_n}}\st\int_{(i-1)\De_n}^{i\De_n}\Big(
\rho_{\si_{s}}(f((i-1)\De_n,X_{s},.)) -
 \rho_{\si_{s}}(f(s,X_{s},.))\Big)\,ds,\\
U'^n_t(3) & = & \frac1{\sqrt{\De_n}}\int_{[t/\De_n]\De_n}^t
\rho_{\si_{s}}(f(s,X_{s},.))\,ds.\end{eqnarray*}
Since  $f$ is  at most with polynomial growth,
$$\frac1{\sqrt{\De_n}}
\int_{[t/\De_n]\De_n}^t\left|\rho_{\si_{s}}(f(s,X_{s},.))\right|\,ds~\le~
K\sqrt{\De_n},$$
hence $U'^{n}_t(3)\,\stackrel{u.c.p.}{\longrightarrow}\,0.$
Otherwise Hypothesis $(M_2)$, and in particular (\ref{ep23}), implies
$$\frac1{\sqrt{\De_n}}\st\int_{(i-1)\De_n}^{i\De_n}|\rho_{\si_{s}}
(f((i-1)\De_n,X_s,.))-
\rho_{\si_{s}}(f(s,X_s,.))|\,ds\le Kt\De_n^{\al-1/2},$$ hence
$U'^{n}_t(2)\,\stackrel{u.c.p.}{\longrightarrow}\,0.$

It remains to show that:

\bee\label{e19} U'^{n}_t(1)\,\stackrel{u.c.p.}{\longrightarrow}\,0.\eee

The function   $(z,x)\mapsto f(\om,s,z,x)$ being $C^1$, so is
the application $ (w,z)  \mapsto \rho_{w}(f(s,z,.)).$
Set $F_{n,i}(\om,w,z):=\rho_w(f(\om,(i-1)\De_n,z,.))$ and
$X''_t=X_t-\int_0^tb_sds$ and $\si''_t=\si_t-\int_0^t\Wb'_sds$.
Then we have  $U'^{n}_t(1)=-\sum_{j=1}^3 U_t^{'n}(1,j)$,  where
\begin{eqnarray*}
U'^n_t(1,1) & = & \frac1{\sqrt{\De_n}}\st\int_{(i-1)\De_n}^{i\De_n}
\left(\int_{(i-1)\De_n}^{s}\left(\widetilde{b}'_u
\frac{ \partial F_{n,i}}{\partial w}(\si_{(i-1)\De_n},X_{(i-1)\De_n})\right.
\right.\\
 & & \left.\left.+\,b'_u
 \frac{ \partial F_{n,i}}{\partial z}(\si_{(i-1)\De_n},X_{(i-1)\De_n})\right)
 \,du\right)\,ds \\
U'^n_t(1,2) & = & \frac1{\sqrt{\De_n}}\st\int_{(i-1)\De_n}^{i\De_n}
\left((\si''_s-\si''_{(i-1)\De_n})
\frac{\partial F_{n,i}}{\partial w}(\si_{(i-1)\De_n},X_{(i-1)\De_n})\right.\\
 & &\left. + (X''_s-X''_{(i-1)\De_n})\frac{ \partial F_{n,i}}{\partial z}
 (\si_{(i-1)\De_n},X_{(i-1)\De_n})\right)
 \,ds\\
U'^n_t(1,3) & = & \frac1{\sqrt{\De_n}}\st\int_{(i-1)\De_n}^{i\De_n}
\Big(
F_{n,i}(\si_{s},X_s)-F_{n,i}(\si_{(i-1)\De_n}, X_{(i-1)\De_n})\\
 & &   -(X_s-X_{(i-1)\De_n}) \frac{ \partial F_{n,i}}{\partial z}
 (\si_{(i-1)\De_n},X_{(i-1)\De_n})\\
  & &  -(\si_{s}-\si_{(i-1)\De_n})
  \frac{ \partial F_{n,i}}{\partial w}(\si_{(i-1)\De_n},X_{(i-1)\De_n})
  \Big)\,ds.
\end{eqnarray*}
Since  $\widetilde{b}',\,b'$  are bounded we have
$\sup_{s\le t} |U'^n_t(1,1)|\,\le  Kt\De_n^{1/2},$
hence $U'^n_t(1,1)\,\stackrel{u.c.p.}{\longrightarrow}\,0.$

Next, the process $U'^n(1,2)$ is a martingale  with respect to the filtration
$(\f_{[t/\De_n]\De_n})$ and the expectation of its  predictable bracket
is smaller than $ Kt\De_n$. Hence   Doob's inequality yields
$U_t^n(1,2)\stackrel{u.c.p.}{\longrightarrow}\,0.$

Finally, if $\ze^n_i(s)$ denotes the integrand in the definition
of $U'^n_t(1,3)$, we have
$$\ze_i^n(s) :=  (\si_{s}-\si_{(i-1)\De_n})\left( \frac{ \partial F_{n,i}}
{\partial w}(\overline{\si}(i,n,s),\overline{X}(i,n,s))-
\frac{ \partial F_{n,i}}{\partial w}(\si_{(i-1)\De_n},X_{(i-1)\De_n}
)\right)$$
   \bee\label{ep46}+\,(X_{s}-X_{(i-1)\De_n})\left( \frac{
   \partial F_{n,i}}{\partial z}(\overline{\si}(i,n,s),\overline{X}(i,n,s))-
   \frac{ \partial F_{n,i}}{\partial z}(\si_{(i-1)\De_n},X_{(i-1)\De_n} )
   \right),\eee
with   $(\overline{\si}(i,n,s),\overline{X}(i,n,s))$ in between
$(\si_{(i-1)\De_n},X_{(i-1)\De_n}) $ and $(\si_{s},X_{s})$.
For $A,\,\ep>0$, set:
$$\begin{array}{l} G_{t}(\ep,A)  =\,   \sup\Big\{\,
\Big|\frac{\partial f}{\partial x}(s,z_1,x_1)-
\frac{\partial f}{\partial x}(s,z_2,x_2)\Big|+
\Big|\frac{\partial f}{\partial z} (s,z_1,x_1)-
\frac{\partial f}{\partial z}(s,z_2,x_2)\Big|: \\[1.8mm]
\qquad\qquad\qquad ~s\le t ;~~|x_1|,~|x_2|\le A;~~|x_1-x_2|\le \ep;~ |z_1|,
~|z_2|\le
K;~~|z_1-z_2|\le \ep\, \Big\},\end{array} $$
then by the properties of $f$, we have $G_t(\ep,A)\to 0$ when $\ep\to 0$.
Therefore it follows from (\ref{ep46}) that
\begin{eqnarray*} |\ze_i^n(s)| & \le & K\left( (1+A)G_t(A\ep,K A)+
\frac{|\si_{s}-\si_{(i-1)\De_n}|+ |X_{s}-X_{(i-1)\De_n}|}{\ep}+ \right. \\
& & \left. + \left(\PP(|U|>A/K)\right)^{1/2}\right)
\ti \Big(|\si_{s}-\si_{(i-1)\De_n}|+ |X_{s}-X_{(i-1)\De_n}|\Big),
\end{eqnarray*} where  $U$ is a  $\n(0,1)$ Gaussian variable.

 Since
under $(LN_2(2))$,  $\E\left\{ |\si_t-\si_s|^2 + |Z_t-Z_s|^2\right\}
\le K|t-s|,$
we  deduce:
\begin{eqnarray*} \frac1{\sqrt{\De_n}}\st\int_{(i-1)\De_n}^{i\De_n}
\E\{|\ze_i^n(s)|\}\,ds
& \le &  Kt\Big(  (1+A) (\E\{G_t(A\ep,KA)^2\})^{1/2}+ \\
 & &  + \PP(|U|>A/K)^{1/2}
+\frac{\sqrt{\De_n}}{\ep}\,\Big).\end{eqnarray*}
Letting $n\to\infty$, then $\ep\to0$, and $A\to\infty$, we obtain
$U'^n_t(1,3)\stackrel{u.c.p.}{\longrightarrow}\,0$, hence (\ref{e19}).\qed
\vspace{.5cm}

The next lemmas are very important
because they deal with the part of the processes having a non-trivial
limit. We use the notation of Subsection \ref{ss1}.
The first one is about the ''even case'' for $f$. Set
\bee\label{ep29}
 \overline{U}^n_t = \sqrt{\De_n}\,\st
 \Big(f((i-1)\De_n,X_{(i-1)\De_n},\be_i^n)
-\E\{f((i-1)\De_n,X_{(i-1)\De_n},\be_i^n)\,|\f_{(i-1)\De_n}\}\,\Big). \eee

\begin{lem}\label{l5}
Suppose   $(LN_2(2))$ satisfied and  $f(\om,s,z,x)$
even in $x$ with at most polynomial growth.
Then $\overline{U}^n_t\stackrel{{\cal L}-(s)}{\longrightarrow} L(f)_t$,
where $L(f)_t=\int_0^t a_s\,d\overline{W}_s$ is given by
(\ref{e35}).
\end{lem}
{\bf Proof:}
Set  $$\xi_i^n=\sqrt{\De_n}\Big(f((i-1)\De_n,X_{(i-1)\De_n},\be_i^n)-
\E\{f((i-1)\De_n,X_{(i-1)\De_n},\be_i^n)|\f_{(i-1)\De_n}\}\Big),$$
then
\bee\label{eol1} E\{\xi_i^n~|\f_{(i-1)\De_n}\}~=~0.\eee
For any   $j,k\in \{1,\cdots,q\}$, we have:

\begin{eqnarray*} \E\{\xi_i^{n,j}\xi_i^{n,k}|\f_{(i-1)\De_n}\} &=&
\De_n\Big(\rho_{\si_{(i-1)\De_n}}((f^jf^k)((i-1)\De_n,X_{(i-1)\De_n},.)
) \\[2mm]
& & -\,\rho_{\si_{(i-1)\De_n}}(f^j((i-1)\De_n,X_{(i-1)\De_n},.))\ti\\[2mm]
& & \ti\rho_{\si_{(i-1)\De_n}}(f^k((i-1)\De_n,X_{(i-1)\De_n},.)))\Big).
\end{eqnarray*}
Then as in  lemma \ref{L44}, one shows that:
\bee\label{eol2} \left.
\begin{array}{l}\st\E\left\{\left(\xi_i^{n,j}\xi_i^{n,k}\right)|
\f_{(i-1)\De_n}\right\}~~\mbox{converges u.c.p. to the process}\\[2mm]
\int_0^t\left(\rho_{\si_{s}}((f^jf^k)(s,X_{s},.))
-\rho_{\si_{s}}(f^j(s,X_{s},.))\rho_{\si_{s}}(f^k(s,X_{s},.))\right)\,ds
\end{array}\right\}\eee

Next for  any $\ep>0$, we have:
\bee\label{eol4} \st \E\{||\xi_i^{n}||^2 1_{\{||\xi_i^n||>\ep\}}
|\f_{(i-1)\De_n}\}~\le~\frac1{\ep^2}\,\st \E\{||\xi_i^{n}||^4
|\f_{(i-1)\De_n}\}~\le~\frac{Kt}{\ep^2}\,\De_n.\eee
Since $f$ is even in $x$:  $\forall\,j'\in\{1,\cdots,m\}$,
\bee\label{eol3} \E\left\{\xi^n_i\De_i^nW^{j'}~|\f_{(i-1)\De_n}\right\}~=~0.
\eee
If now $N$ is a martingale orthogonal to $W$, by the proof of Proposition
4.1 (see (4.13)) of \cite{BGJPN},
\bee\label{ep45} \E\{\xi_i^n\De_i^nN|\f_{(i-1)\De_n}\}=0.\eee

By  (\ref{eol1}), (\ref{eol2}), (\ref{eol4}), (\ref{eol3}) and (\ref{ep45})
we can apply  theorem IX-7-28  of \cite{JJ} which gives our  lemma. \qed

\vspace{.2cm}

\begin{rem}
In the previous lemma, the hypothesis on $f$ is more than what we need,
having $f(\om,s,z,x)$ to be optional even in x and
satisfying $(K(\R^d))$ and  with  at most polynomial growth would
be enough.
\end{rem}

Now we deal with  the case where $f(\om,s,z,x)$ is not even in $x$.

\begin{lem}\label{l9}
Suppose that $X$ and $f$ satisfy respectively  $(LN_2(2))$ and $(LM_2)$,  then
$\overline{U}^n_t\stackrel{{\cal L}-(s)}{\longrightarrow} L(f)_t$, where
$L(f)_t$ is given by (\ref{e20}).
\end{lem}

\ni{\bf Proof:} The proof goes as for lemma \ref{l5}, except
that (\ref{eol3}) fails here, since  $f(\om,s,z,x)$ is not even in $x$.
However we have
$$  \E\Big\{\xi^{n,j}_i\De_i^nW^{k}|\f_{(i-1)\De_n}\Big\}=
\sqrt{\De_n}\E\Big\{f^j((i-1)\De_n,X_{(i-1)\De_n},\be_i^n)
\De_i^nW^{k})|\f_{(i-1)\De_n}\Big\},$$
and (as in the proof of lemma \ref{L44}) one has:
\bee\label{ep24}\st\E\left\{\xi^{n,j}_i\De_i^nW^{k}\,|\f_{(i-1)\De_n}\right\}\,
\stackrel{u.c.p.}{\longrightarrow}
\,\int_0^t w(1)^{j,k}_s\,ds.\eee
Then taking account (\ref{ep24}), and using once more theorem IX-7-28  of
\cite{JJ}, we get this time Lemma \ref{l9}. \qed

\subsubsection{ Proof of theorems \ref{t5} and \ref{t6}:}
We first prove the theorems under the strong hypotheses stated at
the beginning of the Subsection \ref{s2}.
Set $$W^n_t\,:=\,\sqrt{\De_n}\left(V'^{n}_t-\int_0^t\,\rho_{\si_{s}}
(f(s,X_{s},.))\,ds\,\right).$$
Then, using the notation  (\ref{ep43}), (\ref{c4e35}),
(\ref{c4e36}) and (\ref{ep29}), we have:
$$W^n_t\,=\,\sqrt{\De_n}\st(L^n_i-\E\{L^n_i|\f_{(i-1)\De_n}\})\,
+\,\overline{U}^n_t\,+\,U^n_t+U'^n_t.$$

The process $\sqrt{\De_n}\st(L^n_i-\E\{L^ni|\f_{(i-1)\De_n}\})$ is a
martingale
with respect to the filtration $(\f_{[t/\De_n]\De_n})$,
whose predictable bracket is smaller
than $\De_n\E\{\|L^n_i\|^2|\f_{(i-1)\De_n}\}$. Hence Lemma \ref{l8} and
Doob's inequality yield that
$$ \sqrt{\De_n}\st(L^n_i-\E\{L^n_i|\f_{(i-1)\De_n}\})\,
\stackrel{u.c.p.}{\longrightarrow}\,0.$$
Moreover $U^n\,\stackrel{u.c.p.}{\longrightarrow}\,0$
by Lemmas \ref{l4} or \ref{l7},
depending on the case. Next, Lemma \ref{L446} yields
$U'^n\stackrel{u.c.p.}{\longrightarrow}0$.
Finally Lemma \ref{l5} for Theorem \ref{t5} and Lemma \ref{l9} for Theorem
\ref{t6}  give that
$\overline{U}^n$ converges stably in law  to the process $L(f)$
given respectively by (\ref{e35}) and (\ref{e20}).

At this stage, we have proved the theorems under the strong assumptions
mentioned above. The general case is deduced by a ''localization''
procedure.

\end{document}